\documentclass[9pt]{article}
\usepackage{amsmath}
\usepackage[dvips]{graphicx}
\usepackage{textcomp}
\usepackage{amsbsy}
\usepackage{latexsym}
\usepackage[mathscr]{eucal}
\usepackage{amsfonts}
\usepackage{amssymb}
\usepackage[usenames]{color}
\usepackage{stmaryrd}

\usepackage{algorithm,algorithmic}

\usepackage{epsf}

\usepackage{stackrel}

\def\RR{\rm \hbox{I\kern-.2em\hbox{R}}}
\def\NN{\rm \hbox{I\kern-.2em\hbox{N}}}
\def\ZZ{\rm {{\rm Z}\kern-.28em{\rm Z}}}
\def\CC{\rm \hbox{C\kern -.5em {\raise .32ex \hbox{$\scriptscriptstyle
|$}}\kern
-.22em{\raise .6ex \hbox{$\scriptscriptstyle |$}}\kern .4em}}

\def\<{\langle}
\def\>{\rangle}
\def\t{\tilde}

\def\e{\varepsilon}

\def\nl{\newline}
\def\o{\overline}

\def\bv{{\bf v}}
\def\ba{{\bf a}}

\def\bn{{\bf n}}
\def\bb{{\bf b}}
\def\bG{{\bf G}}
\def\bX{{\bf X}}
\def\bI{{\bf I}}
\def\bd{{\bf d}}

\def\cO{{\cal O}}
\newcommand{\R}{{\mathbb R}}
\newcommand{\N}{{\mathbb N}}
\def\Chi{\raise .3ex
\hbox{\large $\chi$}} 
\def\lsima{\hbox{\kern -.6em\raisebox{-1ex}{$~\stackrel{\textstyle<}{\sim}~$}}\kern -.4em}
\def\lsim{\hbox{\kern -.2em\raisebox{-1ex}{$~\stackrel{\textstyle<}{\sim}~$}}\kern -.2em}
\def\[{\Bigl [}
\def\]{\Bigr ]}
\def\({\Bigl (}
\def\){\Bigr )}
\def\[{\Bigl [}
\def\]{\Bigr ]}
\def\({\Bigl (}
\def\){\Bigr )}

\newcommand{\be}{\begin{equation}}
\newcommand{\ee}{\end{equation}}
\newcommand{\beu}{\begin{equation*}}
\newcommand{\eeu}{\end{equation*}}
\newcommand{\bea}{$$ \begin{array}{lll}}
\newcommand{\eea}{\end{array} $$}
\newcommand{\bi}{\begin{itemize}}
\newcommand{\ei}{\end{itemize}}
\newcommand{\iref}[1]{{\rm (\ref{#1})}}
\newtheorem{theorem}{Theorem}
\newtheorem{remark}{Remark}

\newtheorem{corollary}{Corollary}

\def\E{\mathbb E}
\def\P{\mathbb P}

\newcommand{\marg}{{\psi}}
\newcommand{\conditional}{{\varphi}}

\DeclareMathOperator*{\argmin}{argmin}

\newcommand{\conditionalupbound}{\Theta}
\newcommand{\constqupbound}{M}

\newcommand{\trunc}{\tau}

\newcommand{\scalanellafigura}{0.32}

\usepackage{fullpage}

\begin{document}
\bibliographystyle{plain}
\title
{
Optimal weighted least-squares methods \footnote{
This research is supported by Institut Universitaire de France 
and the ERC AdV project BREAD.}
}
\author{
Albert Cohen\thanks{
Sorbonne Universit\'es, UPMC Univ Paris 06, CNRS, UMR 7598, Laboratoire Jacques-Louis Lions, 4, place Jussieu 75005, Paris, France. 
email: cohen@ljll.math.upmc.fr
}
\
and 
\
Giovanni Migliorati\thanks{
Sorbonne Universit\'es, UPMC Univ Paris 06, CNRS, UMR 7598, Laboratoire Jacques-Louis Lions, 4, place Jussieu 75005, Paris, France. 
email: migliorati@ljll.math.upmc.fr
} }

\maketitle
\date{}

\begin{abstract}
We consider the problem of reconstructing an unknown bounded function $u$ defined on a domain $X\subset \mathbb{R}^d$ from noiseless or noisy samples of $u$ at $n$ points $(x^i)_{i=1,\dots,n}$.
We measure the reconstruction error in a norm $L^2(X,d\rho)$ for some given probability measure $d\rho$. Given a linear space $V_m$ with ${\rm dim}(V_m)=m\leq n$, we study in general terms the weighted least-squares approximations from the spaces $V_m$ based on independent random samples. 
It is well known that least-squares approximations can be inaccurate and unstable when $m$ is too close to $n$, even in the noiseless case. 
Recent results from \cite{DH,JNZ} have shown the interest of using weighted least squares for reducing the number $n$ of samples that is needed to achieve an accuracy comparable to that of best approximation in $V_m$, compared to standard least squares as studied in \cite{CDL}.
The contribution of the present paper is twofold. From the theoretical perspective, we establish results in expectation and in probability for weighted least squares in general approximation spaces $V_m$. These results show that for an optimal choice of sampling measure $d\mu$ and weight $w$, which depends on the space $V_m$ and on the measure $d\rho$,  stability and optimal accuracy are achieved under the mild condition that $n$ scales linearly with $m$ up to an additional logarithmic factor. In contrast to \cite{CDL}, the present analysis covers cases where the function $u$ and its approximants from $V_m$ are unbounded, which might occur for instance in the relevant case where $X=\R^d$ and $d\rho$ is the Gaussian measure.
From the numerical perspective, we propose a sampling method which allows one to generate independent and identically distributed samples from the optimal measure $d\mu$.
This method becomes of interest in the multivariate setting where $d\mu$ is generally not of tensor product type. We illustrate this for particular examples of approximation spaces $V_m$ of polynomial type, where the domain $X$ is allowed
to be unbounded and high or even infinite dimensional, motivated by certain applications to parametric and stochastic PDEs.
\end{abstract}

{\bf AMS classification numbers:}   
41A10, 
41A25, 
41A65, 
62E17, 
93E24. 

{\bf Keywords:} 
multivariate approximation, 
weighted least squares, 
error analysis,
convergence rates, 
random matrices,
conditional sampling,
polynomial approximation.

\section{Introduction}

Let $X$ be a Borel set of $\R^d$. We consider the problem of estimating an 
unknown function $u: X \to \R$ from pointwise data $(y^i)_{i=1,\dots,n}$ which are either 
noiseless or noisy observations of $u$ at points $(x^i)_{i=1,\dots,n}$ from $X$. In numerous
applications of interest, some prior information is either established or assumed on the function $u$. 
Such information may take various forms such as:
\begin{itemize}
\item[(i)] regularity properties of $u$, in the sense that it belongs to a given smoothness class;
\item[(ii)] decay or sparsity of the expansion of $u$ in some given basis;
\item[(iii)] approximability of $u$ with some prescribed error by given finite-dimensional spaces.
\end{itemize}
Note that the above are often related to one another and sometimes equivalent,
since many smoothness classes can be characterized by prescribed approximation rates when using 
certain finite-dimensional spaces or truncated expansions in certain bases. 

This paper uses the third type of prior information, taking therefore the view that
$u$ can be ``well approximated'' in some space $V_m$ of functions defined
everywhere on $X$, such that $\dim(V_m)=m$. We work under the 
following mild assumption:
\be
for\; any\; x\in X, \; there\; exists\;  v\in V_m\; such\;  that\; v(x)\neq 0.
\label{assumV}
\ee
This assumption holds, for example, when $V_m$ contains the constant functions.
Typically, the space $V_m$ comes
from a family $(V_j)_{j\geq 1}$ of nested spaces with increasing dimension,
such as algebraic or trigonometric polynomials, or piecewise polynomial functions
on a hierarchy of meshes. 

We are interested in measuring the error in the $L^2(X,d\rho)$ norm
\beu
\|v \|:=\(\int_X |v|^2 d\rho\)^{1/2},
\eeu
where $d\rho$ is a given probability measure on $X$. 
We denote by $\<\cdot,\cdot\>$ the associated inner product.
One typical strategy is to pick the estimate from a finite-dimensional space $V_m$ such that $\dim(V_m)=m$. 
The ideal estimator is given by the $L^2(X,d\rho)$ orthogonal projection 
of $u$
onto $V_m$,
namely
\beu
P_m u:= \argmin_{v\in V_m} \|u-v\|.
\eeu
In general,
this estimator is not computable from a finite number of observations. 
The best approximation error
\beu
e_m(u):= \min_{v\in V_m}\|u-v\|=\|u-P_mu\|,
\eeu
thus serves as a benchmark for a numerical method based 
on a finite sample. In the subsequent analysis, we make significant
use of an arbitrary $L^2(X,d\rho)$ orthonormal basis $\{L_1,\dots,L_m\}$ of the
space $V_m$. 
We also introduce the notation
\beu
e_m(u)_\infty:= \min_{v\in V_m}\|u-v\|_{L^\infty},
\eeu
where $L^\infty$ is meant with respect to $d\rho$, and observe that 
$e_m(u)\leq e_m(u)_\infty$ for any probability measure $d\rho$.

The {\it weighted least-squares} method consists in defining the estimator as
\be
u_W:=\argmin_{v\in V_m} \frac 1 n\sum_{i=1}^n w^i |v(x^i)-y^i|^2,
\label{ls}
\ee
where the weights $w^i> 0$ are given. In the noiseless case $y^i=u(x^i)$, this also writes
\be
\argmin_{v\in V_m}\|u-v\|_n,
\label{dls}
\ee
where the discrete seminorm is defined by 
\be
\|v\|_n:= \left( \frac 1 n\sum_{i=1}^n w^i |v(x^i)|^2 \right)^{1/2}.
\label{lsnorm}
\ee

This seminorm is associated with the 
semi-inner product $\langle \cdot,\cdot \rangle_n$.
If we  
expand the solution to \iref{dls} as
$\sum_{j=1}^m v_j L_j$, the vector $\bv=(v_j)_{j=1,\dots,m}$
is the solution to the normal equations
\be
\bG \bv=\bd,
\label{sys}
\ee
where the matrix $\bG$ has entries $\bG_{j,k}=\<L_j,L_k\>_n$ and where
the data vector $\bd=(d_j)_{j=1,\dots,m}$ is given by $d_j:=\frac 1 n\sum_{i=1}^n w^i y^i L_j(x^i)$.
This system always has at least one solution, which is unique when $\bG$ is nonsingular.
When $\bG$ is singular, we may define $u_W$ as the unique minimal $\ell^2$ norm solution to \iref{sys}.

Note that $\bG$ is nonsingular if and only if  $\|\cdot\|_n$ is a proper norm on the space $V_m$. 
Then, if the data are noisefree
that is, when $y^i=u(x^i)$, we may also write
\beu
u_W=P_m^nu,
\eeu
where $P_m^n$ is the orthogonal projection onto $V_m$ for the norm $\|\cdot\|_n$.

In practice, for the estimator \eqref{ls} to be easily computable, it is important that the functions $L_1,\ldots,L_m$  have explicit expressions that can be evaluated at any point in $X$
so that the system \iref{sys} can be assembled.
Let us note that computing this estimator by solving \iref{sys} 
only requires that $\{L_1,\dots,L_m\}$ is
a basis of the space $V_m$, not necessarily orthonormal in $L^2(X,d\rho)$.
Yet, since our subsequent analysis of this estimator makes
use of an $L^2(X,d\rho)$ orthonormal basis, we simply assume that
$\{L_1,\dots,L_m\}$ is of such type.

In our subsequent analysis, we sometimes work under the 
assumption of a known uniform bound 
\be
\|u\|_{L^\infty}\leq \trunc.
\label{unibound}
\ee
We introduce the truncation operator
\beu
z\mapsto T_{\trunc} (z):={\rm sign}(z)\min\{|z|,\trunc\},
\eeu
and we study the {\it truncated weighted least-squares approximation} defined by 
\beu
u_T:=T_{\trunc}\circ u_W.
\eeu
Note that, in view of \iref{unibound}, we have $|u-u_T|Ê\leq |u-u_W|$ in the pointwise
sense and therefore
\beu
\|u-u_T\|Ê\leq \|u-u_W\|.
\eeu
The truncation operator aims at avoiding unstabilities which may occur when the 
matrix $\bG$ is ill-conditioned. In this paper, we use randomly chosen points $x^i$,
and corresponding weights $w^i=w(x^i)$, distributed in such a way that the resulting random matrix
$\bG$ concentrates towards the identity $\bI$ as $n$ increases.
Therefore, if no $L^\infty$ bound is known, an alternative 
strategy consists in setting to zero the estimator when $\bG$ deviates from
the identity by more than a given value in the spectral norm.
We recall that for $m\times m$
matrices $\bX$, this norm is
defined as $\|\bX\|_2:=\sup_{\|\bv\|_2=1} \|\bX\bv\|_2$.
More precisely, we introduce the
{\it conditioned least-squares approximation}, defined by  
\beu
u_C:=\begin{cases} 
u_W, & \textrm{if } \|\bG-\bI\|_2\leq \frac 1 2, \\
0,   & \textrm{otherwise}.
\end{cases}
\eeu
The choice of $\frac 1 2$ as a threshold for the distance between $\bG$ and $\bI$ in the
spectral norm is related to
our subsequent analysis. However, the value $\frac12$ could be be replaced by any 
real number in $]0,1[$ up to some minor changes in 
the formulation of our results. Note that 
\be
\|\bG-\bI\|_2\leq \frac 1 2 
\implies
{\rm cond}(\bG)\leq 3.
\label{eq:implication_cond}
\ee

It is well known that if $n\geq m$ is too much close to $m$, weighted least-squares methods
may become unstable and inaccurate for most sampling distributions. For example, if $X=[-1,1]$ and
$V_m=\P_{m-1}$ is the space of algebraic polynomials of degree $m-1$, then with $m=n$
the estimator coincides with the Lagrange polynomial interpolation which can be highly unstable
and inaccurate, in particular for equispaced points. The question that we want to address here
in general terms is therefore:
\nl

 {\it Given a space $V_m$ and a measure $d\rho$, how to best choose
the samples $y^i$ and weights $w^i$ in order to ensure that the $L^2(X,d\rho)$
error $\|u-\t u\|$ is comparable to $e_m(u)$, with $n$ being as close as 
possible to $m$, for $\t u\in \{u_W,u_T,u_C\}$ ?}
\nl
 
We address this question in the case where the $x^i$ are
randomly 
chosen. More precisely, we draw independently the $x^i$
according to a certain probabiity measure $d\mu$ defined on $X$.
A natural prescription for the success of the method is that $\|v\|_n$ 
approaches $\|v\|$ as $n$ tends to $+\infty$. Therefore, one first obvious choice is to use
\be
d\mu=d\rho \quad {\rm and} \quad w^i=1, \quad i=1,\ldots,n,
\label{eq:unweighted_ls_w_sigma}
\ee
that is, sample according to the measure in which we plan to evaluate the $L^2$ error
and use equal weights. 
When using equal weights $w^i=1$, the weighted least-squares estimator \eqref{ls} becomes the  \emph{standard least-squares} estimator, as a particular case.
The strategy \eqref{eq:unweighted_ls_w_sigma} was analyzed in \cite{CDL}, through the introduction of the function
\beu
x\mapsto k_m(x):=\sum_{j=1}^m|L_j(x)|^2,
\eeu
which is the diagonal of the integral kernel of the projector $P_m$. This function
only depends on $V_m$ and $d\rho$. It is strictly positive in $X$ due to Assumption \ref{assumV}. Its reciprocal function is characterized by
\beu
\frac {1 }{k_m(x)}=\min_{v\in V_m,v(x)=1} \|v\|^2,
\eeu
and is called 
Christoffel function in the particular case where $V_m$ is the space 
of algebraic polynomials of total degree $m-1$, see \cite{N86}.  Obviously, the function $k_m$ satisfies
\be
\int_X k_m d\rho=m.
\label{intDm}
\ee
We define
\beu
K_m=K_m(V_m,d\rho):=\|k_m\|_{L^\infty},
\eeu
and recall the following results from \cite{CDL,MNST} for the standard least-squares method with the 
weights and the sampling measure chosen as in \eqref{eq:unweighted_ls_w_sigma}.

\begin{theorem}
\label{theo1}
For any $r>0$, if $m$ and $n$ are such that the condition 
\be
K_m \leq \kappa\frac n {\ln n},\;\; {\rm with}\;\; \kappa:=\kappa(r)=\frac{1-\ln 2} {2+2r}
\label{condm}
\ee
is satisfied, then the following hold:
\begin{enumerate}
\item[{\rm (i)}] The matrix $\bG$ satisfies the tail bound
\be
{\rm Pr} \, \left\{ \|\bG-\bI\|_2 > \frac 1 2 \right\} \leq
2n^{-r}. \label{tailhalf}
\ee
\item[{\rm (ii)}] If $u\in L^2(X,d\rho)$ satisfies a uniform bound \iref{unibound}, then
the truncated least-squares estimator satisfies, in the noiseless case,
\be
\E(\|u-u_T\|^2)\leq (1+\e(n))e_m(u)^2+8\trunc^2n^{-r},
\label{nearoptimal}
\ee
where $\e(n):= \frac {4\kappa} {\ln (n)}\to 0$ as $n\to +\infty$, and $\kappa$ as in \iref{condm}.
\item[{\rm (iii)}] If $u\in L^\infty(X,d\rho)$, then the truncated and nontruncated 
least-squares estimators satisfy, in the noiseless case,
\be
\|u- u_T\|\leq \|u-u_W\|\leq (1+\sqrt 2) e_m(u)_\infty,
\label{nearoptimalprob}
\ee
with probability larger than $1-2n^{-r}$.
\end{enumerate}
\end{theorem}

The second item in the above result shows that the optimal accuracy $e_m(u)$
is met in expectation, up to an additional term of order $n^{-r}$. When $e_m(u)$ has polynomial
decay $\cO(m^{-s})$, we are ensured that this additional term can be made negligible
by taking $r$ strictly larger than $s/2$, which amounts in taking $\kappa(r)$ small enough.
Condition \iref{condm} imposes a minimal number of samples to ensure stability and accuracy of standard least squares.
Since  \iref{intDm} implies that $K_m\geq m$, the fulfillment of this
condition requires that $n$ is at least of the order $m\ln(m)$. However simple examples
show that the restriction can be more severe, for example if $V_m=\P_{m-1}$ on $X=[-1,1]$
and with $\rho$ being the uniform probability measure. In this case, 
one choice for the $L_j$ are
the Legendre polynomials with proper normalization $\|L_j\|_{L^\infty}=|L_j(1)|=\sqrt {1+2j}$ so that $K_m=m^2$,
and therefore condition
\iref{condm} imposes that $n$ is at least of order $m^2\ln(m)$.
Other examples in the multivariate setting are discussed in \cite{CCD,CCMNT}
which show that for many relevant approximation spaces $V_m$ and probability
measures $d\rho$, the behaviour of $K_m$ 
is superlinear in $m$, leading to a very demanding regime in terms of the
needed number $n$ of samples. 
In the case of multivariate downward closed polynomial spaces, precise upper bounds for $K_m$ have been proven in \cite{CCMNT,M2015} for measures associated to Jacobi polynomials. 
In addition, note that the above theory does not cover simple
situations such as algebraic polynomials over unbounded domains, for example $X=\R$ equipped with the Gaussian measure, since the 
orthonormal polynomials $L_j$ are unbounded for $j\geq 2$ and thus $K_m=\infty$ if $m\geq 2$.

\section{Main results}

In the present paper, we show that these limitations 
can be overcome, by using a proper
weighted least-squares method. 
We thus return to the general form of the 
discrete norm \iref{lsnorm} used in the definition of the weighted least-squares estimator.
We now use a sampling measure $d\mu$ which generally differs from $d\rho$ and is such that
\beu
w d\mu =  d\rho,
\eeu
where $w$ is a positive function defined everywhere on $X$ and such that $\int_X w^{-1} d\rho=1$, and we then consider the weighted least-square method with weights given by
\beu
w^i=w(x^i).
\eeu
With such a choice, the norm $\|v\|_n$ again approaches $\|v\|$ as $n$ increases. 
The particular case $d\mu=d\rho$ and $w\equiv 1$ corresponds to the standard 
least-squares method analyzed by Theorem \ref{theo1}. Note that changing 
the 
sampling measure is a commonly used strategy for reducing the variance 
in Monte Carlo methods, where it is referred to as {\it importance sampling}.

With $L_j$ again denoting
the $L^2(X,d\rho)$ orthonormal basis of $V_m$, we now introduce the function
\beu
x\mapsto k_{m,w}(x):=\sum_{j=1}^mw(x)|L_j(x)|^2,
\eeu
which only depends on $V_m$, $d\rho$ and $w$,
as well as
\beu
K_{m,w}=K_{m,w}(V_m,d\rho,w):= \|k_{m,w} \|_{L^\infty}.
\eeu
Note that, since the $\sqrt w L_j$ are an $L^2(X,d\mu)$ orthonormal basis of $\sqrt w V_m$,
we find that $\int_X k_{m,w}d\mu=m$ and thus $K_{m,w} \geq m$.
We prove in this paper the following generalization of Theorem \ref{theo1}.

\begin{theorem}
\label{theo2}
For any $r>0$, if $m$ and $n$ are such that the condition 
\be
K_{m,w} \leq \kappa\frac n {\ln n},\;\; {\rm with}\;\; \kappa:=\frac{1-\ln 2} {2+2r}\label{condmw}
\ee
is satisfied, then the following hold:
\begin{enumerate}
\item[{\rm (i)}] The matrix $\bG$ satisfies the tail bound
\be
{\rm Pr} \, \left\{ \|\bG-\bI\|_2 > \frac 1 2 \right\} \leq
2n^{-r}. 
\label{tailhalfw}
\ee
\item[{\rm (ii)}] If $u\in L^2(X,d\rho)$ satisfies a uniform bound \iref{unibound}, then
the truncated weighted least-squares estimator satisfies, in the noiseless case,
\be
\E(\|u-u_T\|^2)\leq (1+\e(n))e_m(u)^2+8\trunc^2n^{-r},
\label{nearoptimalww}
\ee
where $\e(n):= \frac {4\kappa} {\ln (n)}\to 0$ as $n\to +\infty$, and $\kappa$ as in {\rm \iref{condm}}.
\item[{\rm (iii)}] If $u\in L^\infty(X,d\rho)$, then the nontruncated weighted
least-squares estimators satisfy, in the noiseless case,
\be
\|u-u_W\|\leq (1+\sqrt 2) e_m(u)_\infty,
\label{nearoptimalprobw}
\ee
with probability larger than $1-2n^{-r}$.
\item[{\rm (iv)}] If $u\in L^2(X,d\rho)$, then
the conditioned weighted least-squares estimator satisfies, in the noiseless case,
\be
\E(\|u-u_C\|^2)\leq (1+\e(n))e_m(u)^2+2\|u\|^2 n^{-r},
\label{nearoptimalcw}
\ee
where $\e(n):= \frac {4\kappa} {\ln (n)}\to 0$ as $n\to +\infty$, and $\kappa$ as in {\rm \iref{condm}}.
\end{enumerate}
\end{theorem}

Let us mention that the quantity $K_{m,w}$ has been considered in \cite{DH}, 
where similar stability and approximation results have been formulated
in a slightly different form (see in particular Theorem 2.1 therein), 
in the specific framework of total degree polynomial  spaces.

The interest of Theorem \ref{theo2} is that it leads us in a natural way
to an optimal sampling strategy for the weighted least-square method. 
We simply take
\be
w:=\frac {m}{k_m}=\frac m {\sum_{j=1}^m |L_j|^2},
\label{idealw}
\ee
and with such a choice for $w$ one readily checks that 
\be
d\mu:=\frac {k_m} m d\rho,
\label{idealsigma}
\ee
is a probability measure on $X$ since $\int_X k_m d\rho=m$. 

In addition, we have for 
this particular choice that
\beu
k_{m,w}=w k_m=m,
\eeu
and therefore 
\beu
K_{m,w}=m.
\eeu
We thus obtain the following result as a consequence of Theorem 2, which shows that 
the above choice of $w$ and $d\mu$ allows us to obtain near-optimal 
estimates for the truncated weighted least-squares estimator, under the 
minimal condition that $n$ is at least of the order $m\ln(m)$.

\begin{corollary}
\label{cor1}
For any $r>0$, if $m$ and $n$ are such that the condition 
\be
m \leq \kappa\frac n {\ln n},\;\; {\rm with}\;\; \kappa:=\frac{1-\ln 2} {2+2r}
\label{condmwc}
\ee
is satisfied, then the conclusions {\rm(i)}, {\rm(ii)}, {\rm(iii)} and {\rm(iv)} of Theorem \ref{theo2} hold for weighted least squares
with the choice of $w$ and $d\mu$ given by \iref{idealw} and \iref{idealsigma}.
\end{corollary}

One of the interests of the above optimal sampling strategy
is that it applies to polynomial approximation
on unbounded domains that were not covered by Theorem \ref{theo1}, in particular
$X=\R$ equipped with the Gaussian measure. In this case, the relevant target
functions $u$ are often nonuniformly bounded and therefore the results
in items (ii) and (iii) of Theorem \ref{theo2}
do not apply. 
The result in item (iv)
for the conditioned estimator $u_C$ remains valid, since it does not
require uniform boundedness of $u$. 

Let us remark that all the above results are independent of the dimension
$d$ of the domain $X$. 
However,
raising $d$ has the unavoidable effect
of restricting the classes of functions for which the best approximation error $e_m(u)$
or $e_m(u)_\infty$ have some prescribed decay, due to the well-known curse
of dimensionality.

Note that the optimal pair $(d\mu,w)$ described by \iref{idealw} and \iref{idealsigma} 
depends on $V_m$, that is
\beu
w=w_m\quad {\rm and} \quad d\mu=d\mu_m.
\eeu
This raises a difficulty for properly choosing the samples
in settings where the choice of $V_m$ is not fixed a-priori, such as in adaptive methods.
In certain particular cases, it is known that $w_m$ and $d\mu_m$ admit
limits $w^*$ and $d\mu^*$ as $m\to \infty$ and are globally equivalent to 
these limits. One typical example is given by the univariate polynomial
spaces $V_m=\P_{m-1}$, when $X=[-1,1]$ and $d\rho=\rho dx$ where $\rho$ is a Jacobi weight and $dx$ is the Lebesgue measure on $X$. 
In this case $d\mu^*$ is the pluripotential equilibrium measure
\beu
d\mu^*=\frac {dx} {2\pi \sqrt{1-x^2}},
\eeu
see e.g. \cite{Tot,ST}, and one has
\beu
cd\mu^* \leq d\mu_m \leq Cd\mu^*, \quad m\geq 1,
\eeu
for some fixed constants $0<c< C<\infty$. Thus, in such a case, the above corollary 
also holds for the choice $w=w^*$ and $d\mu=d\mu^*$
under the condition $m \leq \frac c C\kappa\frac n {\ln n}$. The development of sampling strategies
in cases of varying values of $m$ without such asymptotic equivalences  
is the object of current investigation.

A closely related weighted least-squares strategy was recently 
proposed and analyzed in \cite{JNZ}, in the polynomial framework. There, the 
authors propose to use the renormalized Christoffel function \iref{idealw} in the definition
of the weights, however sampling from the fixed pluripotential equilibrium measure $d\mu^*$.
Due to the fact that $d\mu_m$ differs from $d\mu^*$, the main estimate obtained in \cite{JNZ} (see p.3 therein) 
does not have the same simple form of a direct comparison between $\|u- u_T \|$ 
and $e_m(u)$ as in {\rm (ii)} of Theorem \ref{theo2}. In particular, it involves an extra 
term $d(f)$ which does not vanish even as $n\to \infty$. 

One intrinsic difficulty when using the optimal pair $(d\mu,w)=(d\mu_m,w_m)$ described by \iref{idealw} and \iref{idealsigma} 
is the effective sample generation, in particular in the multivariate framework since the measure
$d\mu_m$ is generally not of tensor product type. 
One possible approach is to use Markov Chain Monte Carlo methods
such as the Metropolis-Hastings algorithm, as explored in \cite{DH}. In such methods
the samples  are mutually correlated, and only 
asymptotically distributed according to the desired sampling measure.
One contribution of the present paper is to propose a straightforward and 
effective sampling strategy 
for generating an arbitrary finite number $n$ of independent samples identically distributed
 according to $d\mu_m$. This strategy requires that $d\rho$ has
tensor product structure and that the spaces $V_m$ are spanned by tensor product bases, 
such as for multivariate polynomial spaces, in which case $d\mu_m$
is generally not of tensor product type.

The rest of our paper is organized as follows. The proof of Theorem \ref{theo2}
is given in \S 3  in a concise form since it follows the same lines as the original results on
standard least squares from \cite{CDL,MNST}.
We devote \S 4 to analog results in the case of samples affected by additive noise,
proving that the estimates are robust under condition \iref{condmw}.
The proposed method for sampling the optimal measure $d\mu_m$ is discussed in \S 5, and we illustrate its effectiveness in \S 6 
by numerical examples.

\section{Proof of Theorem \ref{theo2}}

The proof is structurally similar to that of Theorem \ref{theo1} given in \cite{CDL} for items (i) and (ii)
and in \cite{CCMNT} for item (iii), therefore we only sketch it. We observe
that $\bG=\frac 1 n \sum_{i=1}^n \bX_i$ where the $\bX_i$ are i.i.d. copies of the 
rank $1$ random matrix
\beu
\bX=\bX(x):=\(w(x)L_j(x)L_k(x)\)_{j,k=1,\dots,m},
\eeu
with $x$ a random variable distributed over $X$ according to $\mu$. One obviously has $\E(\bX)=\bI$.
We then invoke the Chernov bound from \cite{T} to obtain that if $\|\bX\|_2\leq R$ almost surely, then,
for any $0<\delta<1$,
\be
{\rm Pr}\, \{\|\bG-\bI\|_2 > \delta\} \leq 2 m \(\frac
{e^{-\delta}}{(1-\delta)^{1-\delta}}\)^{1/R} =2m \exp\(-\frac
{c_\delta} R\),
\label{eq:finer_est_prob}
\ee
with $c_\delta:=\delta+(1-\delta)\ln(1-\delta)>0$. Taking $\delta=\frac 1 2$, and observing that
\beu
\|\bX(x)\|_2=\frac 1 n w(x)\sum_{j=1}^m|L_j(x)|^2=\frac {{K_{m,w}}(x)} n,
\eeu
we may thus take $R=\frac {K_{m,w}}n$ which yields \iref{tailhalfw} in item (i).

For the proof of \iref{nearoptimalww} in item (ii), we first 
consider the event where $\|\bG-\bI\|_2 \leq \frac 1 2$. In this case we write
\beu
\|u- u_T\|^2=Ê\|T_\trunc (u)-T_\trunc (u_W)\|^2 \leq  \|u-u_W\|^2 =\|u-P_m^n u\|^2 \leq \|g\|^2+\|P_m^n g\|^2, \quad g:=u-P_mu,
\eeu
where we have used that $P_m^nP_m u=P_m u$ and that $g$ is orthogonal to $V_m$,
and thus
\beu
\|u-u_T\|^2  \leq e_m(u)^2+\sum_{j=1}^m |a_j|^2,
\eeu
where $\ba=(a_j)_{j=1,\dots,m}$ is solution of the system
\beu
\bG\ba=\bb, \;\;
\eeu
and $\bb:=(\<g,L_k\>_n)_{k=1,\dots,m}$. Since 
$\|\bG^{-1}\|_2\leq 2$, it follows that
\beu
\|u-u_T\|^2  \leq e_m(u)^2+4 \sum_{k=1}^m |\<g,L_k\>_n|^2.
\eeu
In the event where $\|\bG-\bI\|_2 > \frac 1 2$, we simply write $\|u-u_T\|\leq 2\trunc$. It follows that
\beu
\E(\|u-u_T\|^2)\leq e_m(u)^2+4 \sum_{k=1}^m\E( |\<g,L_k\>_n|^2)+ 8\trunc^2 n^{-r}.
\eeu
For the second term, we have
\begin{align*}
\E( |\<g,L_k\>_n|^2) & =\frac 1 {n^2}\sum_{i=1}^n\sum_{j=1}^n \E(w(x^i)w (x^j)g(x^i)g(x^j)L_k(x^i)L_k(x^j)) \\
 & =\frac 1 {n^2} \(n(n-1)|\E(w(x)g(x)L_k(x))|^2+n\E(|w(x)g(x)L_k(x)|^2)\) \\
 & = \(1-\frac 1 n\)|\<g,L_k\>|^2 +\frac 1 n\int_X|w(x)|^2 |g(x)|^2|L_k(x)|^2d\mu \\
 &=\frac 1 n\int_X w(x)|g(x)|^2|L_k(x)|^2d\rho,
 \end{align*}
where we have used the fact that $g$ is $L^2(X,\rho)$-orthogonal to $V_m$ and thus to $L_k$.
Summing over $k$, we obtain
\beu
\sum_{k=1}^m\E( |\<g,L_k\>_n|^2)\leq \frac {K_{m,w}} n \|g\|^2
\leq \frac \kappa {\ln (n)}e_m(u)^2,
\eeu
and we therefore obtain \iref{nearoptimalww}. 

For the proof of \iref{nearoptimalprobw} in item (iii) we place ourselves in the event
where $\|\bG-\bI\|_2 \leq \frac 1 2$. This property also means that
\beu
\frac 1 2 \|\bv\|^2_2\leq \<\bG \bv, \bv\>_2 \leq \frac 3 2 \|\bv\|^2_2, \quad \bv\in\R^m,
\eeu
which can be expressed as a norm equivalence over $V_m$,
\be
\frac 1 2 \|v\|^2 \leq \|v\|_n^2 \leq \frac 3 2 \|v\|^2, \quad v\in V_m.
\label{normeq}
\ee
We then write that for any $v\in V_m$,
\begin{align*}
\|u-\P_m^n u\| & \leq \|u-v\|+\|v-P_m^n u\| \\
& \leq \|u-v\|+\sqrt 2\|v-P_m^n u\|_n\\
& \leq \|u-v\|+\sqrt 2\|u-v\|_n\\
& \leq (1+\sqrt 2)\|u-v\|_{L^\infty},\\
\end{align*}
where we have used \iref{normeq}, the Pythagorean identity $\|u-v\|_n^2=\|u-P_m^nu\|_n^2+\|v-P_m^n u\|_n^2$,
and the fact that both $\|\cdot\|$ and $\|\cdot\|_n$ are dominated by $\|\cdot\|_{L^\infty}$. Since $v$ is arbitrary,
we obtain \iref{nearoptimalprobw}.

Finally, \iref{nearoptimalcw} in item (iv) is proven in a very similar way as \iref{nearoptimalww} in item (ii),
by writing that in the event $\|\bG-\bI\|_2 > \frac 1 2$, we have $\|u-u_C\|=\|u\|$,
so that
\beu
\E(\|u-u_C\|^2)\leq e_m(u)^2+4 \sum_{k=1}^m\E( |\<g,L_k\>_n|^2)+ 2\|u\|^2 n^{-r},
\eeu
and we conclude in the same way.
\hfill $\Box$

\section{The noisy case}

In a similar way as in \cite{CDL,MNT2015}, we can analyze the case where the observations of $u$
are affected by an additive noise.
In practical situations the noise may come from different sources, such as 
a discretization error when $u$ is evaluated by some numerical code,
or a measurement error.
The first one may be viewed as 
a perturbation of $u$ by a deterministic funtion $h$, that is, we observe
\beu
y^i=u(x^i)+h(x^i).
\eeu
The second one is typically modelled as
a stochastic fluctuation, that is, we observe
\beu
y^i=u(x^i)+\eta^i.
\eeu
where $\eta^i$ are independent realizations of the centered random variable $\eta=y-u(x)$. Here, we do not
necessarily assume $\eta$ and $x$ to be independent, however we typically
assume that the noise is centered, that is,
\be
\E(\eta| x)=0,
\ee
and we also assume uniformly bounded conditional variance
\be
\label{condMeanVariance}
\sigma^2:=\sup_{x \in X} \E (|\eta|^2 |x)<\infty.
\ee
Note that we may also consider consider a noncentered noise, which amounts
in adding the two contributions, that is,
\be
y^i=u(x^i)+\beta^i, \quad \beta^i=h(x^i)+\eta^i,
\label{noise3}
\ee
with $h(x)=\E(\beta|x)$. The following result shows that the estimates in Theorem \ref{theo2}
are robust under the presence of such an additive noise.

\begin{theorem}
\label{theo3}
For any $r>0$, if $m$ and $n$ are such that condition \iref{condmw}
is satisfied, then the following hold for the noise model \iref{noise3}:
\begin{enumerate}
\item[{\rm (i)}] if $u\in L^2(X,d\rho)$ satisfies a uniform bound \iref{unibound}, then
the truncated weighted least-squares estimator satisfies
\be
\E(\|u-u_T\|^2)\leq (1+ 2\e(n))e_m(u)^2+(8+ 2\e(n))\|h\|^2+\frac {\o K_{m,w}\sigma^2}{n}+8\trunc^2n^{-r},
\label{nearoptimalwwnoise1}
\ee
\item[{\rm (ii)}] if $u\in L^2(X,d\rho)$, then
the conditioned weighted least-squares estimator satisfies
\be
\E(\|u-u_C\|^2)\leq  (1+ 2\e(n))e_m(u)^2+(8+ 2\e(n))\|h\|^2+\frac {\o K_{m,w}\sigma^2}{n}+2\|u\|^2 n^{-r},
\label{nearoptimalcwnosie1}
\ee
\end{enumerate}
where in both cases $\e(n):= \frac {4\kappa} {\ln (n)}\to 0$ as $n\to +\infty$, 
with $\kappa$ as in {\rm \iref{condm}}, and
$\o K_{m,w}:=\int_X k_{m,w}d\rho$.
\end{theorem}
\noindent
{\bf Proof:} We again first consider the event where $\|\bG-\bI\|_2 \leq \frac 1 2$. In this case we write
\beu
\|u- u_T\| \leq \|u-u_W\|, 
\eeu
and use the decomposition
$u-u_W=g-P^n_m g- h$
where $g=u+P_mu$ as in the proof of Theorem \ref{theo2} and 
$h$ stands for the solution to the least-squares problem
for the noise data $(\beta^i)_{i=1,\dots,n}$. Therefore
$$
\|u-u_W\|^2 =\|g\|^2+\|P_m^n g+h\|^2 \leq \|g\|^2+2\|P_m^n
g\|^2+2\|h\|^2 =\|g\|^2+2\|P_m^n
g\|^2 +2\sum_{j=1}^m
|n_j|^2,
$$
where $\bn=(n_j)_{j=1,\dots,m}$ is
solution to 
$$
\bG\bn=\bb, \quad \bb:=\(\frac 1 n\sum_{i=1}^n\beta^i w(x^i)L_k(x^i)\)_{k=1,\dots,m}=(b_k)_{k=1,\dots,m}.
$$
Since
$\|\bG^{-1}\|_2\leq 2$, it follows that
\beu
\|u-u_T\|^2  \leq e_m(u)^2+ 8\sum_{k=1}^m |\<g,L_k\>_n|^2+8\sum_{k=1}^m|b_k|^2.
\eeu
Compared to the proof of Theorem \ref{theo2}, we need to estimate the
expectation of the third term on the right side. For this we simply write that
$$
\E(|b_k|^2)=\frac 1 {n^2}\sum_{i=1}^n\sum_{j=1}^n \E(\beta^i w(x^i)L_k(x^i)\beta^j w(x^j)L_k(x^j)).
$$
For $i\neq j$, we have
$$
\E(\beta^iw(x^i)L_k(x^i)\beta^j w(x^j)L_k(x^j))=\E(\beta w(x)L_k(x))^2=\E(h(x) w(x)L_k(x))^2
=\left |\int_X h wL_k d\mu\right |^2=|\<h,L_k\>|^2.
$$
Note that the first and second expectations are with respect to the joint density of $(x,\beta)$ and
the third one with respect to the density of $x$, that is, $\mu$.
For $i=j$, we have
\begin{align*}
\E(|\beta^i w(x^i)L_k(x^i)|^2)&=\E(|\beta w(x) L_k(x)|^2) \\
&=\int_X \E(|\beta w(x) L_k(x)|^2 | x) d\mu\\
&=\int_X\E(|\beta|^2 |x)|w(x)L_k(x)|^2d\mu\\
&=\int_X\E(|\beta|^2 |x) w(x)|L_k(x)|^2d\rho\\
&=\int_X(|h(x)|^2+\E(|\eta|^2Ê|x)) w(x)|L_k(x)|^2d\rho\\
&\leq  \int_X(|h(x)|^2 +\sigma^2)w(x)|L_k(x)|^2d\rho.
\end{align*}
Summing up on $i$, $j$ and $k$, and using condition \iref{condmw},
we obtain that 
\be
\sum_{k=1}^m \E(|b_k|^2)\leq \(1-\frac 1{n^2}\) \|h\|^2+ \frac{K_{m,w}}{n}\|h\|^2+
\frac {\o K_{m,w}}{n}\sigma^2 
\leq 
\(1+\frac{\kappa}{\log n}\) \|h\|^2+ \frac {\o K_{m,w}\sigma^2}{n}.
\ee
For the rest we proceed as for item (ii) and (iv) in the proof of Theorem \ref{theo2}, using 
that in the event $\|\bG-\bI\|_2 > \frac 1 2$ we have 
$\|u-u_T\|\leq 2\trunc$ and $\|u-u_C\|=\|u\|$.
\hfill $\Box$

\begin{remark}
Note that for the standard least-squares method, corresponding to the case where $w\equiv 1$, we know that
$\o K_{m,w}=m$. The noise term thus
takes the stardard form $ \frac {m\sigma^2}n$, as seen for example in
Theorem~3 of \cite{CDL} or in Theorem~1 of \cite{MNT2015}. Note that, in any case,
condition  \iref{condmw} implies that this term is 
bounded by $\frac{\kappa\sigma^2}{\log n}$.
\end{remark}

The conclusions of Theorem \ref{theo3} do not include the 
estimate in probability similar to item (iii) in Theorem \ref{theo2}. 
We can obtain such an estimate in
the case of a bounded noise, where we assume that $h\in L^\infty(X)$ and $\eta$
is a bounded random variable, or equivalently, assuming that $\beta$ is a bounded
random variable, that is we use the noise model \iref{noise3} with
\be
|\beta| \leq D, \quad a.s.
\label{noiseuni}
\ee
For this bounded noise model we have the following result.

\begin{theorem}
\label{theo4}
For any $r>0$, if $m$ and $n$ are such that condition \iref{condmw}
is satisfied, then the following hold
for the the noise model \iref{noise3} under \iref{noiseuni}:
if $u\in L^\infty(X,d\rho)$, then the nontruncated weighted
least-squares estimator satisfies
\be
\|u-u_W\|\leq (1+\sqrt 2) e_m(u)_\infty+ \sqrt 2 D,
\label{nearoptimalprobwnoise}
\ee
with probability larger than $1-2n^{-r}$.
\end{theorem}
\noindent
{\bf Proof:} Similar to the proof of (iii) in Theorem \ref{theo2}, we place ourselves in the event
where $\|\bG-\bI\|_2 \leq \frac 1 2$ and use the norm equivalence 
\iref{normeq}.
We then write that for any $v\in V_m$,
$$
\|u-u_W\| \leq \|u-v\|+\|v-P_m^n u \|+\| P_m^n \beta\|.
$$
The first two terms already appeared in the noiseless
case and can be treated in the same way. The new term $P_m^n\beta$ corresponds to the weighted least-squares approximation from
the noise vector, and satisfies
$$
\| P_m^n \beta\| \leq \sqrt 2 \| P_m^n \beta\|_n \leq \sqrt 2 \|\beta\|_n\leq \sqrt 2 D.
$$
This leads to \iref{nearoptimalprobwnoise}.
\hfill $\Box$

\section{Random sampling from $\mu_m$}

The analysis in the previous sections prescribes the use of the optimal sampling measure $d\mu_m$ defined in \eqref{idealsigma} for drawing the samples 
$x^1,\ldots,x^n$ in the weighted least-squares method. 
In this section we discuss numerical methods for generating independent random samples according to this measure, in a specific relevant multivariate setting.

Here, we make the assumption that $X=\times_{i=1}^d X_i$ is a Cartesian product of univariate real domains $X_i$, and that $d\rho$ is a product measure, that is,
\beu
d\rho=\bigotimes_{i=1}^d d\rho_i,
\eeu
where
each $d\rho_i$ is a measure defined on $X_i$. We assume that each $d\rho_i$ is of the form
\beu
d\rho_i(t)=\rho_i(t) dt,
\eeu
for some nonnegative continuous function $\rho_i$, and therefore
\beu
d\rho(x)= \rho(x)\,dx , \quad \rho(x)= \prod_{i=1}^d \rho_i(x_i), \quad x=(x_1,\dots,x_d)\in X.
\eeu
In particular $d\rho$ is absolutely continuous with respect to the Lebesgue measure.

We consider the following general setting: for each $i=1,\dots,d$, we choose a univariate basis $(\phi_j^i)_{j\geq 0}$ orthonormal in $L^2(X_i,d\rho_i)$. We then define the tensorized basis
\beu
L_\nu(x):=\prod_{i=1}^d \phi_{\nu_i}^i(x_i), \quad \nu\in\N_0^d,
\eeu
which is orthonormal in $L^2(X,d\rho)$. We consider general subspaces of the form
\beu
V_m:={\rm span}\{L_\nu\; : \; \nu\in \Lambda\},
\eeu
for some multi-index set $\Lambda\subset \N_0^d$ such that $\#(\Lambda)=m$. 
Thus we may rename the $(L_\nu)_{\nu\in\Lambda}$ as $(L_j)_{j=1,\dots,m}$ after a proper ordering has been chosen, for example in the lexicographical sense. 
For the given set $\Lambda$ of interest, we introduce 
\beu
\lambda_j:=\max_{\nu \in \Lambda} \nu_j \quad {\rm and} \quad \lambda_\Lambda:=\max_{j=1,\ldots,d} \lambda_j.
\eeu

The measure $d\mu_m$ is thus given by $d\mu_m(x)=\mu_m(x) dx$, where
\begin{equation}
\mu_m(x) := \dfrac{1}{m} \sum_{i=1}^m |L_i(x)|^2 \rho(x)= \dfrac{1}{\#(\Lambda)} \sum_{\nu\in\Lambda} |L_\nu(x)|^2 \rho(x), \quad x \in X.
\label{eq:multiv_density_sigma}
\end{equation}

We now discuss our sampling method for generating $n$ independent random samples $x^1,\dots, x^n$ identically distributed according to the multivariate density \eqref{eq:multiv_density_sigma}. 
Note that this density does not have a product structure, despite $\rho$ is a product density. 
There exist many methods for sampling from multivariate densities. 
In contrast to Markov Chain Monte Carlo methods mentioned in the introduction, the method that we next propose exploits the particular structure of the multivariate density \eqref{eq:multiv_density_sigma}, in order to generate independent samples in a straightforward manner, and sampling only from univariate densities.

Given the vector $x=(x_1,\ldots,x_d)$ of all the coordinates, for any $A\subseteq\{1,\ldots,d\}$, we introduce the notation 
\beu
x_{A}:=(x_i)_{i\in A},
\quad 
\bar A:=\{1,\ldots, d\}\setminus A, 
\quad
x_{\bar A}:=(x_i)_{i \in \bar A},
\eeu
and 
\beu
dx_A:=\stackrel[i\in A]{}{\bigotimes} d x_i, \quad d\rho_A :=\stackrel[i\in A]{}{\bigotimes} d\rho_i, \quad \rho_A(x_A):=\prod_{i\in A}Ê\rho_i(x_i), \quad X_A:=\stackrel[i\in A]{}{\times} X_i.
\eeu
In the following, we mainly use the particular sets
\beu
A^q:=\{1,\dots,q\}\quad {\rm and} \quad \bar A^q:=\{q+1,\dots,d\},
\eeu
so that any $x\in X $ may be  written as $x=(x_{A^q},x_{\bar A^q})$. 

Using such a notation, for any $q=1,\ldots,d$, we associate to the joint density $\mu_m$ its marginal density $\marg_q$ of the first $q$ variables, namely 
\begin{equation}
\marg_{q}(x_{A^q}):=\int_{X_{\bar A^q}} \mu_m(x_{A^q},x_{\bar A^q})\, dx_{\bar A^q}. 
\label{eq:def_marginal_densities}
\end{equation}
Since $(\phi_j^i)_{j\geq 0}$ is an orthonormal basis of $L^2(X_i,d\rho_i)$, for any $q=1,\ldots,d$ and any $\nu\in\N_0^d$, we obtain that 
\beu
\int_{X_{\bar A^q}} |L_\nu (x_{A^q},x_{\bar A^q})|^2 \rho(x_{A^q},x_{\bar A^q}) dx_{\bar A^q}=\rho_{A^q}(x_{A^q})\prod_{i=1}^{q} | \phi_{\nu_i}^i(x_i) |^2, \quad x_{A^q}\in X_{A^q}.
\eeu
Therefore, the marginal density \eqref{eq:def_marginal_densities} can be written in simple form as 
\be
\marg_{q}(x_{A^q})= \dfrac{1}{\#(\Lambda)}
\rho_{A^q}(x_{A^q})
\sum_{\nu \in \Lambda} 
\prod_{i=1}^{q} | \phi_{\nu_i}^i(x_i) |^2. 
\label{eq:expression_marginals_first_q-1}
\ee

\paragraph{Sequential conditional sampling.}
Based on the previous notation and remarks, we propose an algorithm which generates $n$ samples $x^k=(x_1^k,\ldots,x_d^k)\in X$ with $k=1,\ldots,n$, that are independent and identically distributed realizations from the density $\mu_m$ in \eqref{eq:multiv_density_sigma}. 

In the multivariate case the coordinates can be arbitrarily reordered. 
Start with the first coordinate $x_1$ and sample $n$ points $x_1^1,\ldots,x_1^n \in X_1$ from the univariate density 
\begin{equation}
\conditional_1:X_1\to\R:t\mapsto
\conditional_1(t):=\marg_1(t)= \dfrac{\rho_1(t)}{\#(\Lambda)} 
\sum_{\nu \in \Lambda} | \phi_{\nu_1}^1(t) |^2,
\label{eq:conditional_x1}
\end{equation}
which coincides with the marginal $\marg_1$ of $x_1$ calculated in 
\eqref{eq:expression_marginals_first_q-1}. 
In the univariate case $d=1$ the algorithm terminates. In the multivariate case $d\geq 2$, by iterating $q$ from $2$ to $d$, consider the $q$th coordinate $x_q$, and sample $n$ points $x_q^1,\ldots,x_q^n \in X_q$ in the following way: 
for any $k=1,\ldots,n$, given the values $x_{A^{q-1}}^k=(x_1^k,\ldots,x_{q-1}^k) \in X_{A^{q-1}}$ that have been 
calculated at the previous $q-1$ steps, sample the point $x_q^k \in X_q $ from the univariate density 
\begin{align}
\conditional_{q}:X_q \to \R:
t\mapsto
\conditional_{q}
(t | 
x_{A^{q-1}}^k 
) 
:=  
\rho_q(t)
\dfrac{
\sum_{\nu \in \Lambda} 
| \phi_{\nu_q}^q(t) |^2  
\prod_{j=1}^{q-1} | \phi_{\nu_j}^j(x_j^k) |^2
}{ 
\sum_{\nu \in \Lambda} 
\prod_{j=1}^{q-1} | \phi_{\nu_j}^j(x_j^k) |^2
}
.
\label{eq:conditional_q_simplified}
\end{align}
The expression on 
the right-hand side of 
\eqref{eq:conditional_q_simplified} is continuous at any $t\in X_q$ and at any $x_{A^{q-1}}^k \in X_{A^{q-1}}$. 
Assumption~\ref{assumV} ensures that the denominator of \eqref{eq:conditional_q_simplified} is strictly positive for any possible choice of $x_{A^{q-1}}^k=(x_1^k,\ldots,x_{q-1}^k) \in X_{A^{q-1}}$, 
and also 
ensures that the marginal $\marg_{q-1}$ is strictly positive at any point $x_{A^{q-1}}^k \in X_{A^{q-1}}$ such that $\rho_{A^{q-1}}(x_{A^{q-1}}^k)\neq 0$. 
For any $t \in X_q$ and any $x_{A^{q-1}}^k \in  X_{A^{q-1}}$ such that $\rho_{A^{q-1}}(x_{A^{q-1}}^k)\neq 0$, the density $\conditional_q$ satisfies 
\begin{align}
\conditional_{q}
(t | 
x_{A^{q-1}}^k 
) 
= 
\dfrac{ 
\marg_{q}(x_{A^{q-1}}^k,t)  
}{
\marg_{q-1}(x_{A^{q-1}}^k)  
},
\label{eq:conditional_q}
\end{align}
where the densities $\marg_{q}$ and $\marg_{q-1}$ are the marginals defined in \eqref{eq:def_marginal_densities} and evaluated  at the points $(x_{A^{q-1}}^k,t) \in X_{A^{q}}$ and $x_{A^{q-1}}^k \in X_{A^{q-1}}$, respectively. 
From \eqref{eq:conditional_q}, using \eqref{eq:expression_marginals_first_q-1} and simplifying the term $\rho_{A^{q-1}}(x_{A^{q-1}}^k)=\prod_{j=1}^{q-1}\rho_j(x_j^k)\neq 0$, one obtains the right-hand side of \eqref{eq:conditional_q_simplified}. 
The right-hand side of equation \eqref{eq:conditional_q} is well defined for any $t \in X_q$ and any $x_{A^{q-1}}^k \in  X_{A^{q-1}}$ such that $\rho_{A^{q-1}}(x_{A^{q-1}}^k)\neq 0$, and it is not defined at the points $x_{A^{q-1}}^k \in  X_{A^{q-1}}$ such that $\rho_{A^{q-1}}(x_{A^{q-1}}^k) = 0$ where $\marg_{q-1}(x_{A^{q-1}}^k)$ vanishes. 
Nonetheless, \eqref{eq:conditional_q} has finite limits at any point $(x_{A^{q-1}}^k,t)\in X_{A^q}$, and these limits equal expression \eqref{eq:conditional_q_simplified}.

According to technical terminology, the right-hand side of equation \eqref{eq:conditional_q} is the conditional density of $x_q$ given $x_1,\ldots,x_{q-1}$ with respect to the density $\marg_q$, and $\conditional_q$ is the continuous extension to $X_{A^q}$ of this conditional density.  

The densities $\conditional_1,\ldots,\conditional_d$ defined in \eqref{eq:conditional_x1}--\eqref{eq:conditional_q_simplified} can be concisely rewritten 
for any $q=1,\ldots,d$ as 
\begin{equation}
\conditional_q(t|
x_{A^{q-1}}^k)= \rho_q(t) \sum_{\nu \in \Lambda} \alpha_{\nu}(x_{A^{q-1}}^k) | \phi_{\nu_q}^q(t) |^2, 
\label{eq:conditional_concise}
\end{equation}
where the nonnegative weights $(\alpha_\nu)_{\nu \in \Lambda}$ are defined as 
\beu
\alpha_{\nu}=
\alpha_{\nu}(z_{A^{q-1}})
:=\begin{cases}
\dfrac{1}{\#(\Lambda)}, & \textrm{ if } q=1, \\
\dfrac{ \prod_{j=1}^{q-1} |\phi_{\nu_j}^j(z_j)|^2 }{
\sum_{\nu \in \Lambda} \prod_{j=1}^{q-1} |\phi_{\nu_j}^j(z_j)|^2 
}, & \textrm{ if } 2\leq q \leq d,
\end{cases}
\eeu
for any $z_{A^{q-1}}=(z_1,\ldots,z_{q-1})\in X_{A^{q-1}}$.  
Since 
$
\sum_{\nu \in \Lambda} \alpha_\nu = 1,
$
each density $\conditional_q$ in \eqref{eq:conditional_concise} is a convex combination of the densities $\rho_q |\phi_1^q|^2,\ldots,\rho_q |\phi_{\lambda_q}^q|^2$. 
Note that if the orthonormal basis $(\phi_j^q)_{j \geq 0}$ 
have explicit expressions and can be evaluated at any point in $X_q$, then the same holds for the univariate densities \eqref{eq:conditional_concise}.
In particular, in the polynomial case, for standards univariate densities $\rho_i$ such as uniform, Chebyshev or Gaussian, the orthonormal polynomials $(\phi_j^i)_{j\geq 1}$ have expressions which are explicitely computable, for example by recursion formulas. 

In Algorithm~\ref{algo:sampling} we summarize our sampling method, that sequentially samples the univariate densities \eqref{eq:conditional_concise} to 
generate independent samples from the multivariate density \eqref{eq:multiv_density_sigma}. In the univariate case $d=1$ the algorithm does not run the innermost loop, and only samples from $\conditional_1$.
In the multivariate case $d\geq 2$ the algorithm runs also the innermost loop, and conditionally samples also from $\conditional_2,\ldots,\conditional_d$. 
Our algorithm therefore relies on accurate sampling methods for the relevant univariate densities
\eqref{eq:conditional_concise}. 

\begin{figure}[!htbp]
\centering
\begin{algorithm}[H]
\caption{\small Sequential conditional sampling for $\mu_m$.}
\begin{algorithmic}
\REQUIRE $n$, $d$, $\Lambda$, 
$\rho_i$, $(\phi_j^i)_{j\geq 0}$ for $i=1,\ldots,d$.
\ENSURE $x^1,\ldots,x^n\stackrel[]{\textrm{i.i.d.}}{\sim} \mu_m$.
\FOR{$k=1$ to $n$}
\STATE $\alpha_\nu\gets(\#(\Lambda))^{-1}$, for any $\nu \in \Lambda$.
\vspace{0.1cm}
\STATE Sample $x_1^k$ from $t\mapsto\conditional_1(t)= \rho_1(t) \stackrel[\nu \in \Lambda]{}{\sum} 
\alpha_{\nu}
\,
 |\phi_{\nu_1}^1(t)|^2$.
\FOR{$q=2$ to $d$}
\STATE $\alpha_\nu
\gets
\dfrac{
\stackrel[j=1]{q-1}{\prod} |\phi_{\nu_j}^j(x_j^k)|^2
}{
\stackrel[ \nu \in \Lambda]{}{\sum} 
\
\stackrel[j=1]{q-1}{\prod} 
|\phi_{\nu_j}^j(x_j^k)|^2
}
$, for any $\nu \in \Lambda$. 
\vspace{0.1cm}
\STATE Sample $x_q^k$ from $t\mapsto \conditional_q(t)= \rho_q(t) 
\stackrel[\nu \in \Lambda]{}{\sum} \alpha_{\nu} \, | \phi_{\nu_q}^q(t) |^2$.
\ENDFOR
\STATE $x^k \gets (x_1^k,\ldots,x_d^k)$. 
\ENDFOR
\end{algorithmic}
\label{algo:sampling}
\end{algorithm}
\end{figure}

We close this section by discussing two possible methods for sampling from such densities: 
 \emph{rejection sampling} 
and 
\emph{inversion transform sampling}.  
Both methods equally apply to any univariate density $\conditional_q$, 
and therefore we present them for 
any $q$ arbitrarily chosen from $1$ to $d$.

\paragraph{Rejection sampling (RS).}
For applying this method, one needs to find a suitable univariate density $\conditionalupbound_q$, whose support contains the support of $\conditional_q$, 
 and a suitable real $\constqupbound_q>1$ such that 
\beu
\conditional_q(t) \leq \constqupbound_q \conditionalupbound_q(t), \quad  t \in \textrm{supp}(\conditional_q). 
\eeu
The density $\conditionalupbound_q$ should be easier to sample than $\conditional_q$, \emph{i.e.}~efficient pseudorandom number generators for sampling from $\conditionalupbound_q$ are available.
The value of $\constqupbound_q$ should be the smallest possible. 
For sampling one point from $\conditional_q$ using RS: 
sample a point $z$ from $\conditionalupbound_q$, and sample $u$ from the standard uniform $\mathcal{U}(0,1)$. 
Then check if $u< \conditional_q(z) / \constqupbound_q \conditionalupbound_q(z)$:  
if this is the case then 
accept $z$ as a realization from $\conditional_q$, 
otherwise reject $z$ and restart sampling $z$ and $u$ from beginning. 
On average, acceptance occurs once every $\constqupbound_q$ trials. 
Therefore, for a given $q$, sampling one point from $\conditional_q$
by RS 
 requires on average $\constqupbound_q$ evaluations of the function 
\beu
t \mapsto
\dfrac {  \conditional_q(t) }{ \constqupbound_q \conditionalupbound_q(t)} = 
\dfrac{ \rho_q(t) }{\constqupbound_q \conditionalupbound_q(t) }
\sum_{\nu \in \Lambda} \alpha_{\nu} |\phi_{\nu_q}^q(t)|^2. 
\eeu
This amounts in evaluating $\constqupbound_q$ times the terms $\phi_0^q,\phi_{\lambda_q}^q$ and a subset of the terms $\phi_1^q,\ldots,\phi_{\lambda_q -1}^q$, depending on $\Lambda$. 
The coefficients $\alpha_\nu$ depend on the terms   
$\phi_0^j,\ldots,\phi_{\lambda_j}^j$ for $j=1,\ldots,q-1$, which have been already evaluated when sampling the previous coordinates $1,\ldots,q-1$.
Thus, if we use RS for sampling the univariate densities, 
the overall computational cost of Algorithm~\ref{algo:sampling}
for sampling $n$ points $x^1,\ldots,x^n\in X$ is on average proportional to $n\sum_{q=1}^d \constqupbound_q (\lambda_q+1)$. 

When the basis functions $(\phi_j^q)_{j\geq 0}$ form a bounded orthonormal system, 
an immediate and simple choice of the parameters in the algorithm is 
\begin{equation}
M_q=\max_{\nu \in \Lambda} \Vert \phi_{\nu_q}^q \Vert_{L^\infty}^2, 
\quad \textrm{ and } 
\quad
\conditionalupbound_q(t)=\rho_q(t).
\label{eq:simple_choice_rejection_sampling}
\end{equation} 
With such a choice, we can quantify more precisely the average computational cost of sampling $n$ points in dimension $d$.    
When $(\phi_j^q)_{j\geq 0}$ are the Chebyshev polynomials, whose $L^\infty$ norms satisfy $\Vert \phi_j^q \Vert_{L^\infty}\leq \sqrt2$, we obtain 
the bound $2n\sum_{q=1}^d (\lambda_q +1) \leq 2 n d (\lambda_\Lambda+1) \leq 2ndm$.
When $(\phi_j^q)_{j\geq 0}$ are the Legendre polynomials, whose $L^\infty$ norms satisfy $\Vert \phi_j^q \Vert_{L^\infty} \leq \sqrt{2j+1}$, we have the crude estimate $2n\sum_{q=1}^d (\lambda_q +1)^2 \leq 2 n d (\lambda_\Lambda+1)^2 \leq 2ndm^2$.  
In general, when $(\phi_j^q)_{j\geq 0}$ are Jacobi polynomials, similar upper bounds can be derived, and the dependence of these bounds on $n$ and $d$ is linear.

\paragraph{Inversion transform sampling (ITS).}
Let $\Phi_q:X_q\to [0,1]$ be the cumulative distribution function associated to the univariate density $\conditional_q$. 
In the following, only when using the ITS method, we make the further assumption that $\rho_q$ vanishes at most a finite number of times in $X_q$. Such an assumption is fulfilled in many relevant situations, \emph{e.g.} 
when  $\rho_q$ is the density associated to Jacobi or Hermite polynomials orthonormal in $L^2(X_q,d\rho_q)$. Together with Assumption~\ref{assumV}, this ensures that the function $t\mapsto\Phi_q(t)$ is continuous and strictly increasing on $X_q$.    
Hence $\Phi_q$ is a bijection between $X_q$ and $[0,1]$, and it has a unique inverse $\Phi_q^{-1}:[0,1]\to X_q$ which is continuous and strictly increasing on $[0,1]$. 
Sampling from $\conditional_q$ using ITS can therefore be performed as follows: sample $n$ independent realizations $u^1,\ldots,u^n$ identically distributed according to the standard uniform $\mathcal{U}(0,1)$, and obtain the $n$ independent samples from $\conditional_q$ as $(\Phi_q^{-1}(u^1),\ldots,\Phi_q^{-1}(u^n))$.  

For any $u\in[0,1]$, computing $z=\Phi_q^{-1}(u)\in X_q$ is equivalent to find the unique solution $z\in X_q$ to $\Phi_q(z)=u$. 
This can be executed by elementary root-finding numerical methods, \emph{e.g.} the bisection method or Newton's method. 
In alternative to using root-finding methods, one can build an interpolant operator $\mathcal{I}_q$ of $\Phi_q^{-1}$, and then approximate $\Phi_q^{-1}(u)\approx \mathcal{I}_q(u)$ for any $u\in[0,1]$.
Such an interpolant $\mathcal{I}_q$ can be constructed for example by piecewise linear interpolation, from the data $(\Phi_q(t_1^q),t_1^q), \ldots , (\Phi_q(t_{s_q}^q),t_{s_q}^q)$ at $s_q$ suitable points $t_1^q<\ldots<t_{s_q}^q$ in $X_q$. 

Both root-finding methods and the interpolation method require evaluating the function $\Phi_q$ pointwise in $X_q$.  
In general these evaluations can be computed using standard univariate quadrature formulas. When $(\phi_j^q)_{j\geq 0}$ are orthogonal polynomials, the explicit expression of the primitive of $\conditional_q$ 
can be used for directly evaluating the function $\Phi_q$. 

Finally we discuss the overall computational cost of Algorithm~\ref{algo:sampling} for sampling $n$ points $x^1,\ldots,x^n\in X$ when using ITS for sampling the univariate densities. 
With the bisection method, this overall cost amounts to $n\sum_{q=1}^d \gamma_q W_q$, where $\gamma_q$ is the maximum number of iterations for locating the zero in $X_q$ up to some desired tolerance, and $W_q$ is the computational cost of each  iteration. 
With the interpolation of $\Phi_q^{-1}$, the overall cost amounts to $n$ evaluations of each interpolant $\mathcal{I}_q$, in addition to the cost for building the interpolants which does not depend on $n$.

\section{Examples and numerical illustrations}
\noindent
This section presents the numerical performances of the weighted least-squares method 
compared to the standard least-squares method, in three relevant situations where $d\rho$ can be either the uniform measure, the Chebyshev measure, or the Gaussian measure. 
In each one of these three cases, we choose $w$ and $d\mu$ in the weighted least-squares method from \iref{idealw} and \iref{idealsigma}, as prescribed by our analysis in Corollary~\ref{cor1}.
For standard least squares we choose $w$ and $d\mu$ as in \eqref{eq:unweighted_ls_w_sigma}. 
Our tests focus on the condition number of the Gramian matrix, that quantifies the stability of the linear system \eqref{sys} and the stability of the weighted and standard least-squares estimators.    
A meaningful 
quantity is therefore 
the probability  
\begin{equation}
\Pr\{ \textrm{cond}(\bG)\leq 3\},
\label{eq:numerical_test_cond_prob}
\end{equation}
where, through \eqref{eq:implication_cond}, the value three of the threshold is related to the parameter $\delta=1/2$ in the previous analysis. 
For any $n$ and $m$, from \eqref{eq:implication_cond} the probability \eqref{eq:numerical_test_cond_prob} is larger than $\Pr\{ \| \bG - \bI \|_2 \leq \frac12 \}$. 
From Corollary~\ref{cor1}, under condition \eqref{condmwc} between $n$, $m$ and $r$, the Gramian matrix of  weighted least squares satisfies \eqref{tailhalfw}
and therefore the probability \eqref{eq:numerical_test_cond_prob} is larger than $1-2n^{-r}$.  
For standard least squares, from Theorem~\ref{theo1} the Gramian matrix satisfies \eqref{eq:numerical_test_cond_prob} with probability larger than $1-2n^{-r}$, but under condition \eqref{condm}. 

In the numerical tests the probability \eqref{eq:numerical_test_cond_prob} is approximated by empirical probability, obtained by counting how many times the event $\textrm{cond}(\bG)\leq 3$ occurs when repeating the random sampling one hundred times.

All the examples presented in this section confine to multivariate approximation spaces of polynomial type.  One natural assumption in this case is to require that the set $\Lambda$ is {\it downward closed}, 
that is, satisfies
\beu
\nu\in \Lambda\quad {\rm and} \quad \t \nu\leq \nu \implies \t \nu\in \Lambda,
\eeu
where $\t \nu\leq \nu$ means that $\t \nu_j\leq \nu_j$ for all $i=1,\dots,d$. Then $V_m$
is the polynomial space spanned by the monomials 
\beu
z\mapsto z^\nu:=\prod_{j=1}^d z_j^{\nu_j},
\eeu
and the orthonormal basis $L_\nu$ is provided by taking each $(\phi_j^i)_{j\geq 0}$ to be
a sequence of univariate orthonormal polynomials of $L^2(X_i,d\rho_i)$.

In both the univariate and multivariate forthcoming examples, the random samples from the measure $d\mu_m$ are generated using Algorithm~\ref{algo:sampling}. 
The univariate densities $\conditional_1,\ldots,\conditional_d$ are sampled using the inversion transform sampling method. 
The inverse of the cumulative distribution function is approximated using the  interpolation technique.

\subsection{Univariate examples}
\begin{figure}[!htb]
\begin{tabular}{ccc}
\small $d\rho$ uniform measure  & \small  $d\rho$ Gaussian measure & \small $d\rho$ Chebyshev measure \\
\includegraphics[clip,width=\scalanellafigura\textwidth,viewport=0 0 1080 962]{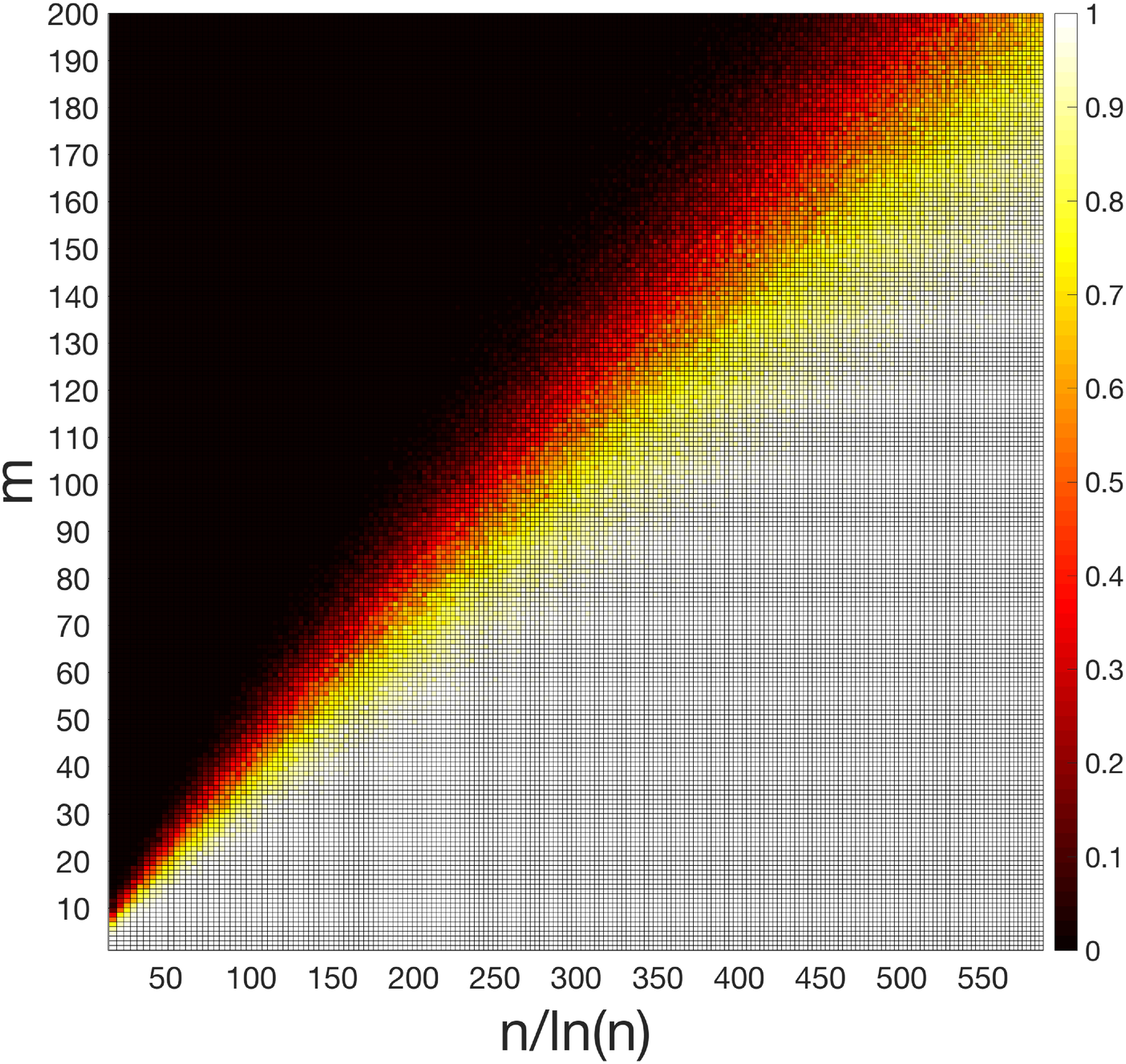}  &
\includegraphics[clip,width=\scalanellafigura\textwidth,viewport=0 0 1080 962]{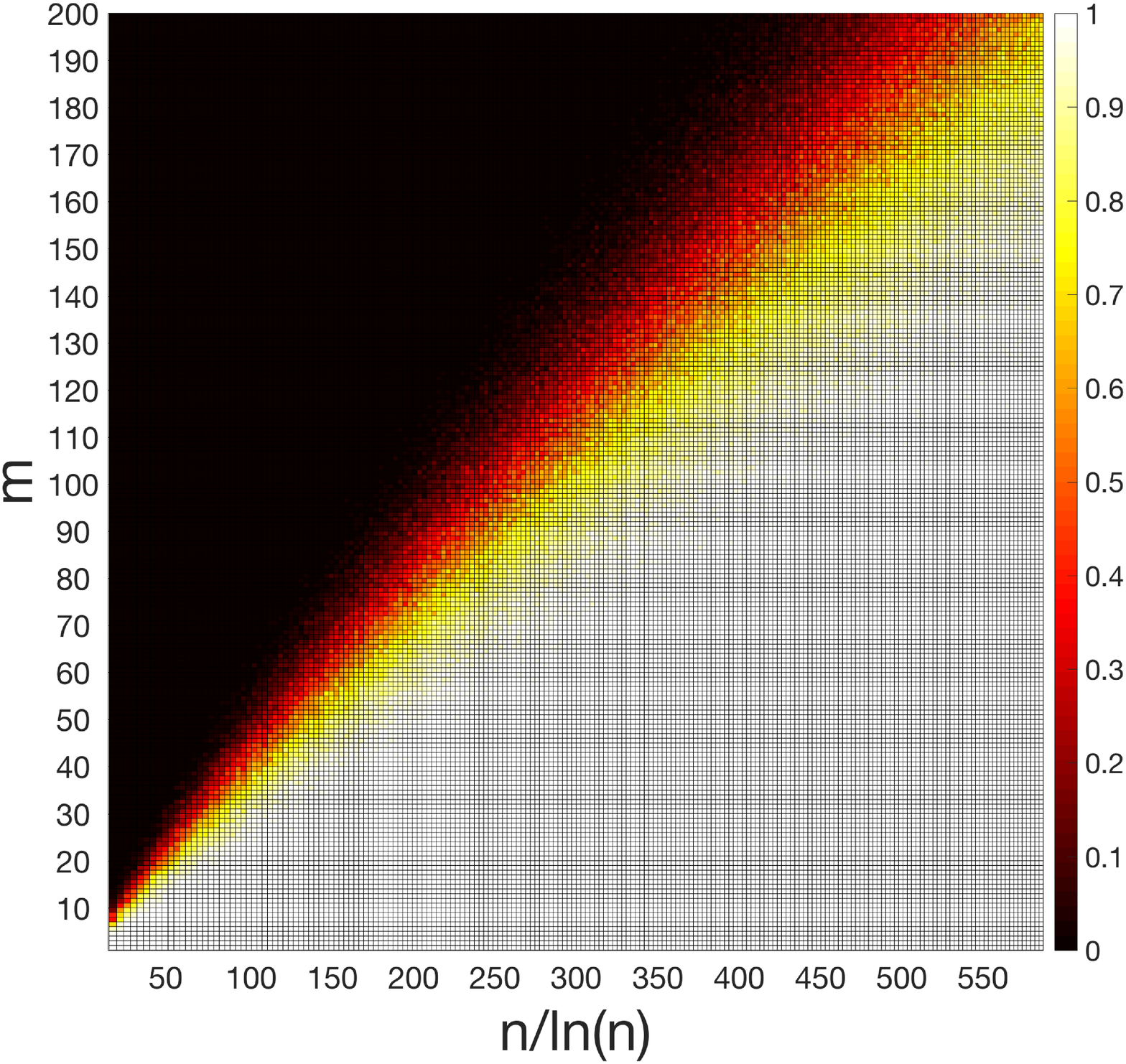}  &
\includegraphics[clip,width=\scalanellafigura\textwidth,viewport=0 0 1080 962]{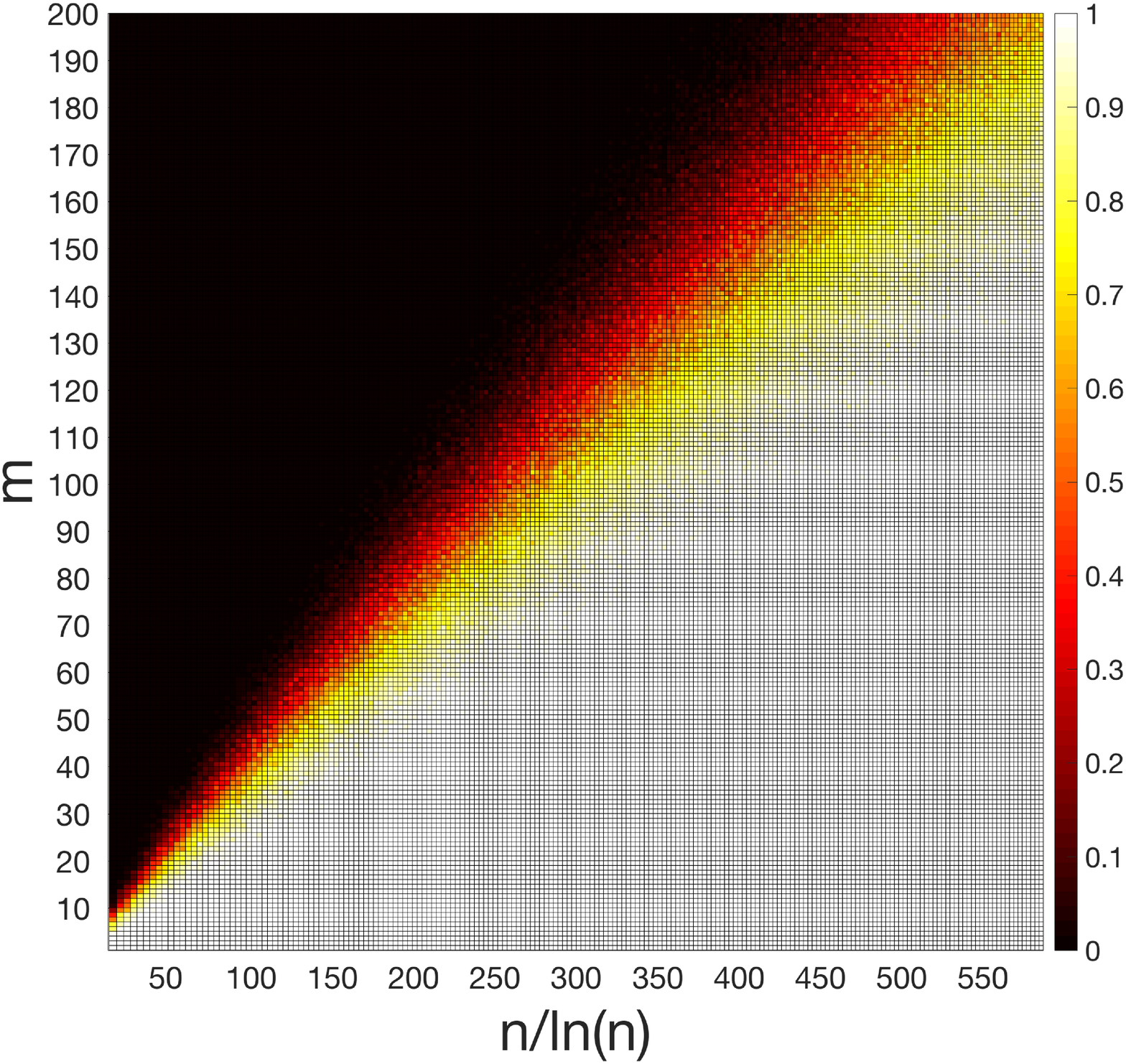}  
\end{tabular}
\caption{
Weighted least squares, 
$Pr\{\textrm{cond}(G)\leq 3\}$, 
$d=1$. 
Left: $d\rho$ uniform measure. 
Center: $d\rho$ Gaussian measure.
Right: $d\rho$ Chebyshev measure. 
}
\label{fig_results_univariate_wls}
\end{figure}

In the univariate case $d=1$, let the index set be $\Lambda=\{0,\ldots,m-1\}$ and $V_m=\P_\Lambda=\textrm{span}\{z^k : k=0,\ldots,m-1 \}$.   
We report in Fig.~\ref{fig_results_univariate_wls} the probability \eqref{eq:numerical_test_cond_prob}, approximated by empirical probability, when $\bG$ is the Gramian matrix of the weighted least-squares method. 
Different combinations of values for $m$ and $n$ are tested, with three choices of the measure $d\rho$: uniform, Gaussian and Chebyshev. 
The results do not show perceivable differences among the performances of weighted least squares with the three different measures. In any of the three cases, $n/\ln(n) \geq 4 m$ is enough to obtain an empirical probability equal to one that $\textrm{cond}(\bG)\leq 3$. This confirms that condition \eqref{condmwc} with any choice of $r>0$ ensures \eqref{eq:numerical_test_cond_prob}, since 
it demands for a larger number of samples.  

\begin{figure}[!htb]
\begin{tabular}{ccc}
\small $d\rho$ uniform measure  & \small  $d\rho$ Gaussian measure & \small $d\rho$ Chebyshev measure \\
\includegraphics[clip,width=\scalanellafigura\textwidth,viewport=0 0 1080 962]{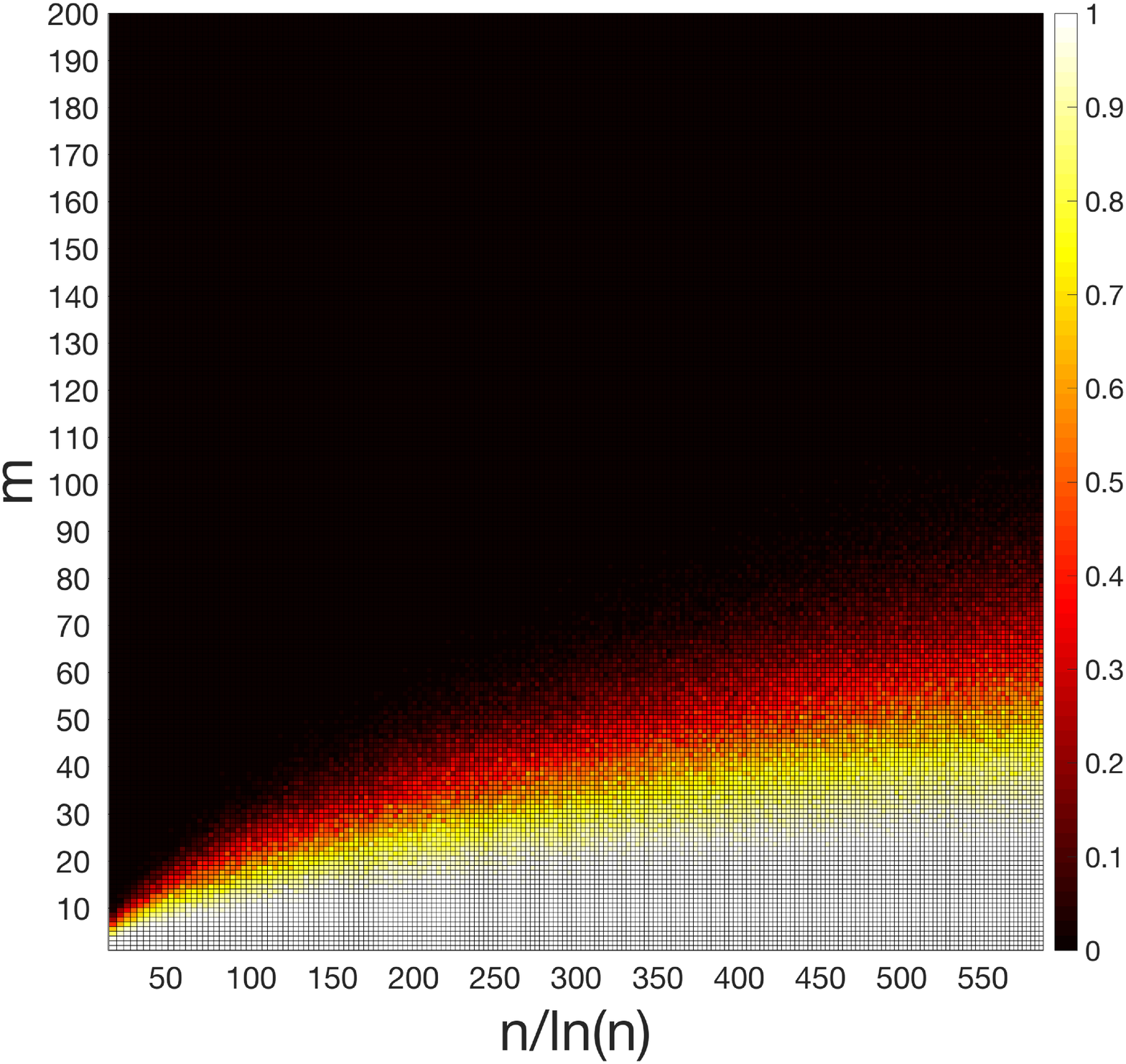}  &
\includegraphics[clip,width=\scalanellafigura\textwidth,viewport=0 0 1080 962]{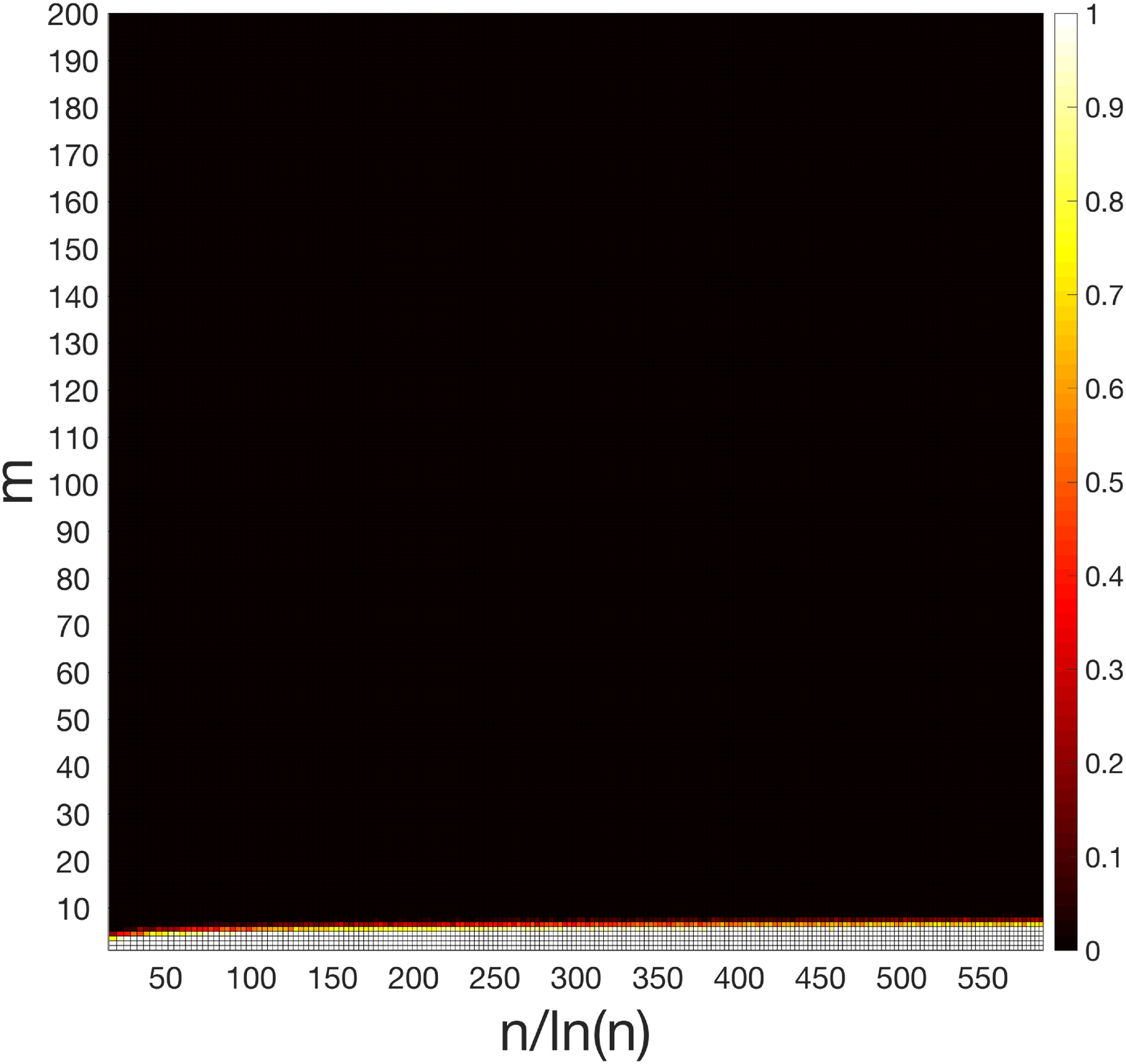}  &
\includegraphics[clip,width=\scalanellafigura\textwidth,viewport=0 0 1080 962]{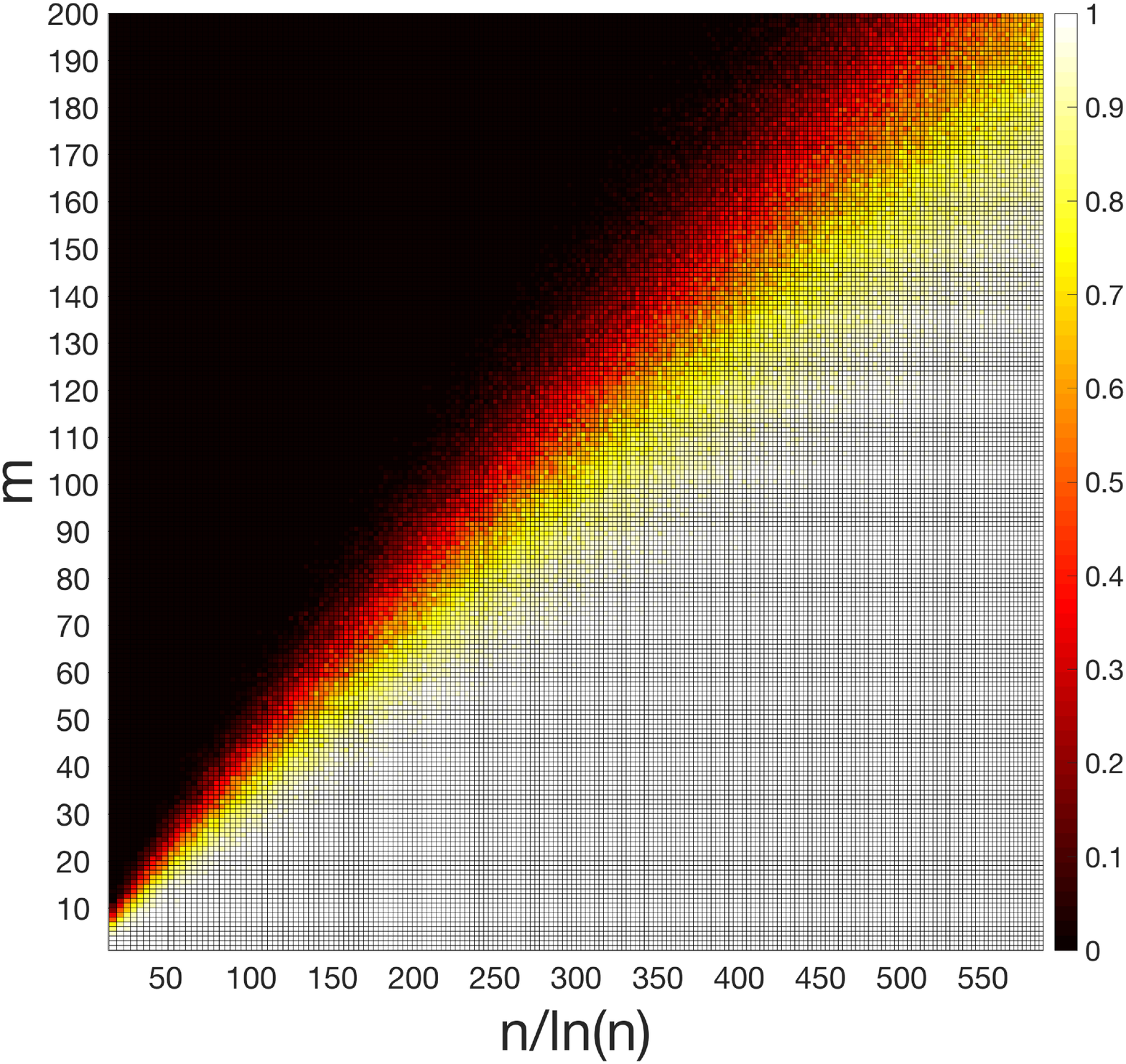}  
\end{tabular}
\caption{
Standard least squares, 
$Pr\{\textrm{cond}(G)\leq 3\}$, 
$d=1$. 
Left: $d\rho$ uniform measure. 
Center: $d\rho$ Gaussian measure.
Right: $d\rho$ Chebyshev measure. 
}
\label{fig_results_univariate_sls}
\end{figure}

Fig.~\ref{fig_results_univariate_sls} shows the probability 
\eqref{eq:numerical_test_cond_prob} when $\bG$ is the Gramian matrix of 
standard least squares. With the uniform measure, the condition $n/\ln(n)\geq m^2$ is enough 
to have \eqref{eq:numerical_test_cond_prob} with empirical probability larger than $0.95$. 
When $d\rho$ is the Gaussian measure, stability requires a very large number of evaluations, roughly $n/\ln(n)$ linearly proportional to $\exp(m/3)$.  
For the univariate Chebyshev measure, 
it is proven that standard least squares are stable under the same minimal condition \eqref{condmwc} as for weighted least squares. 
In accordance with the theory, the numerical results obtained in this case with weighted and standard least squares are indistinguishable, see 
Fig.~\ref{fig_results_univariate_wls}-right and Fig.~\ref{fig_results_univariate_sls}-right.

\subsection{Multivariate examples}

\begin{figure}[!htb]
\begin{tabular}{ccc}
\small $d\rho$ uniform measure  & \small  $d\rho$ Gaussian measure & \small $d\rho$ Chebyshev measure \\
\includegraphics[clip,width=\scalanellafigura\textwidth,viewport=0 0 1080 962]{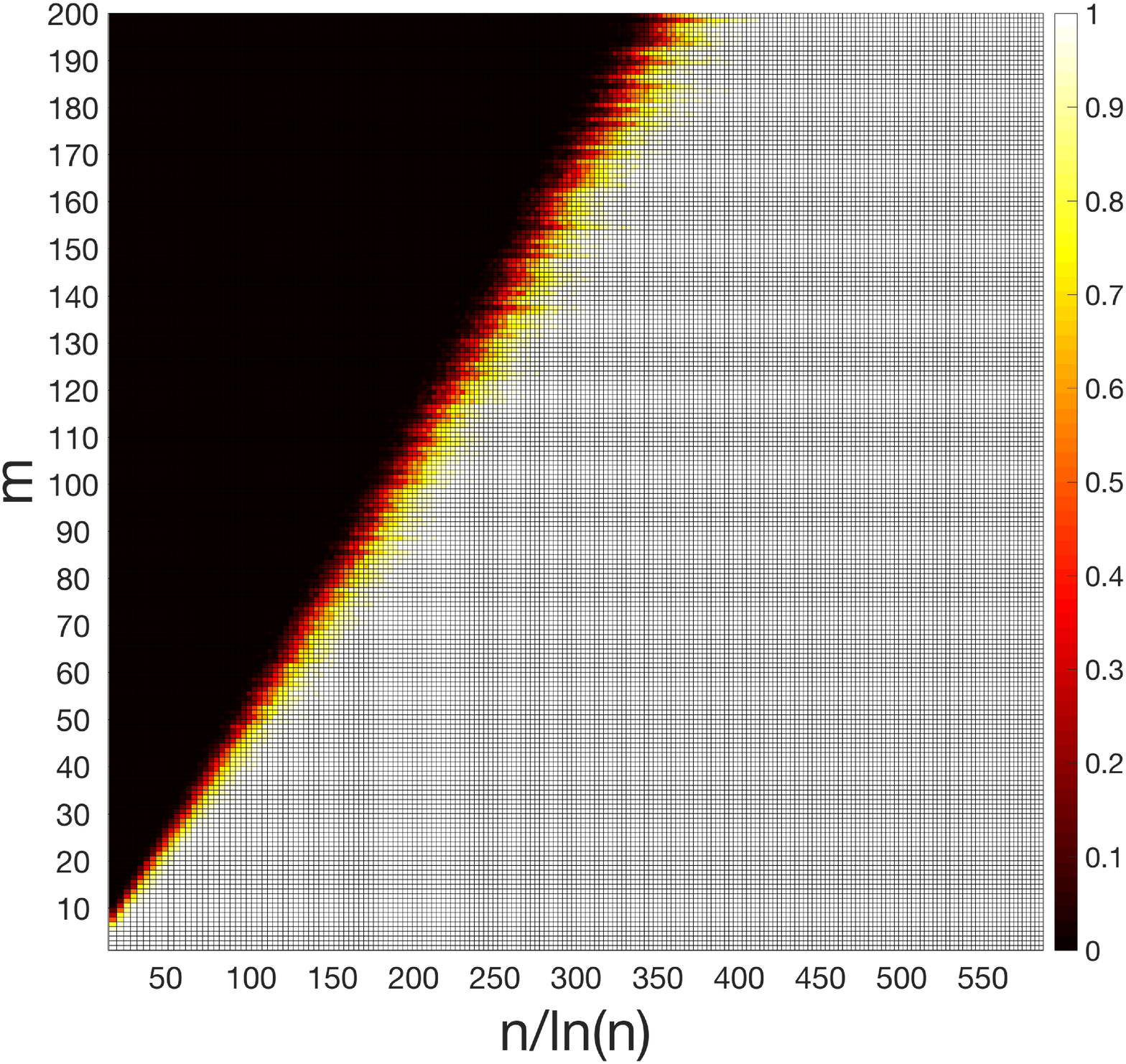} &
\includegraphics[clip,width=\scalanellafigura\textwidth,viewport=0 0 1080 962]{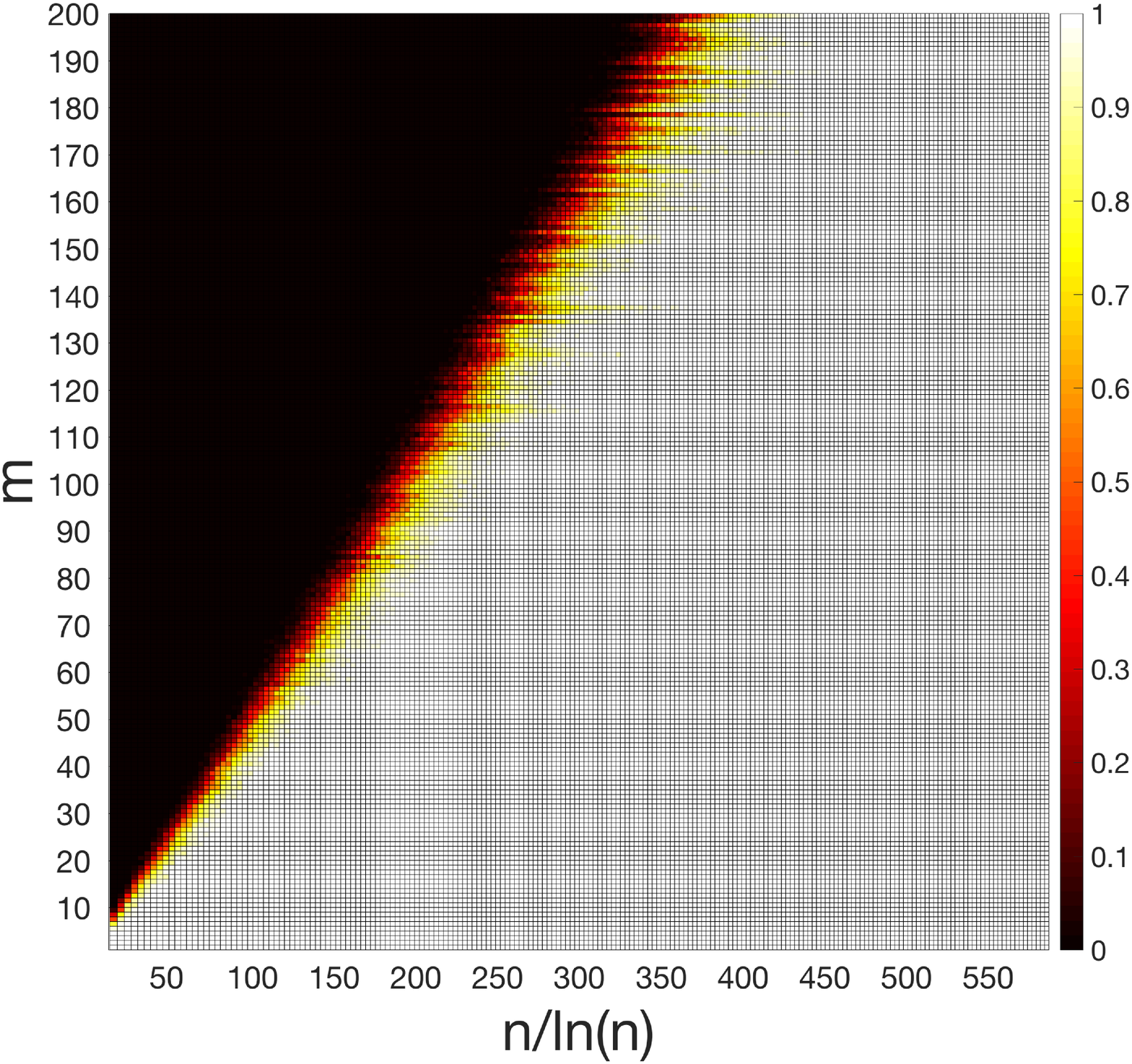}  &
\includegraphics[clip,width=\scalanellafigura\textwidth,viewport=0 0 1080 962]{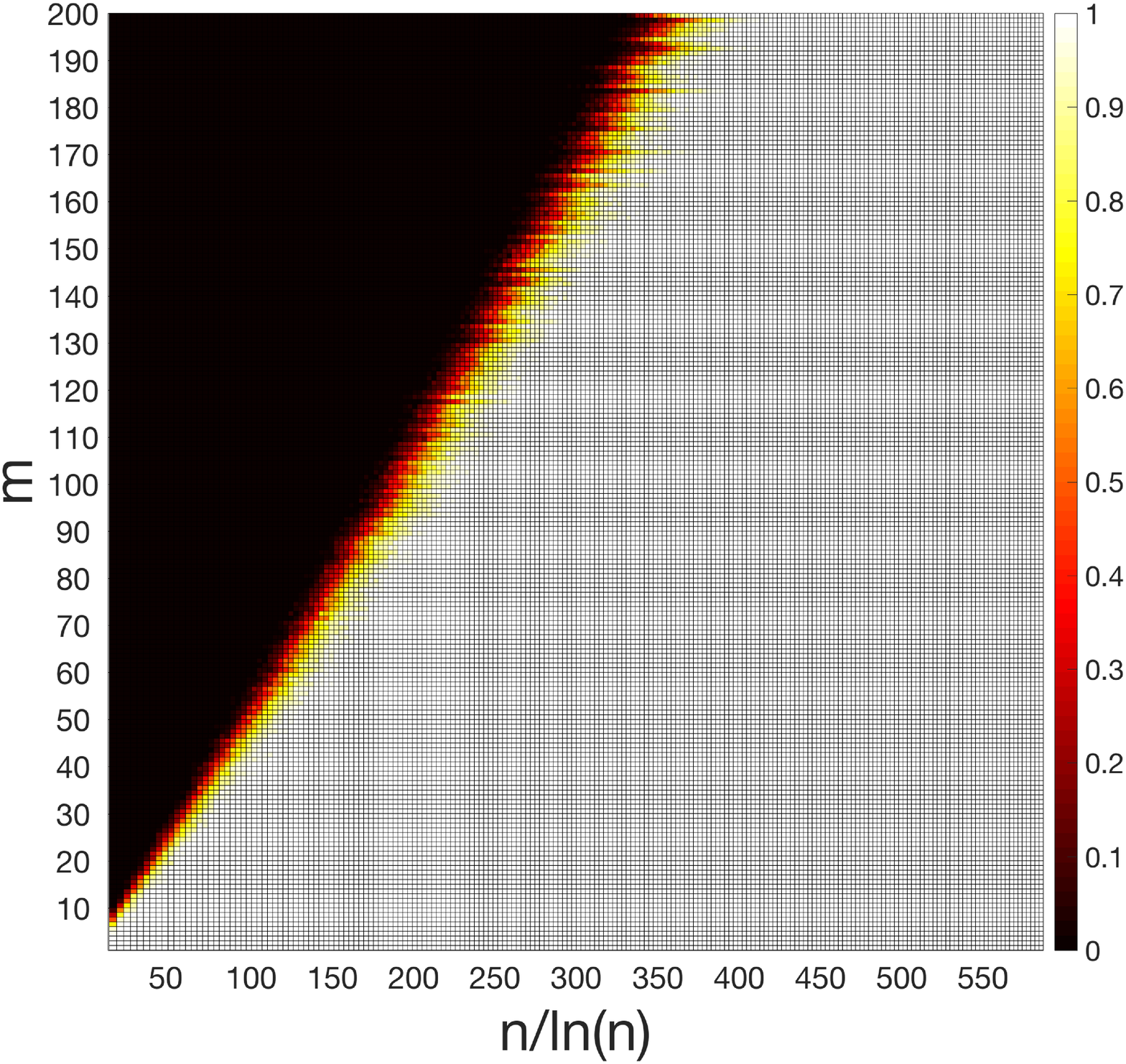}  
\end{tabular}
\caption{
Weighted least squares, 
$Pr\{\textrm{cond}(G)\leq 3\}$, 
$d=10$. 
Left: $d\rho$ uniform measure. 
Center: $d\rho$ Gaussian measure.
Right: $d\rho$ Chebyshev measure. 
}
\label{fig_results_multivariate_d10_wls}
\end{figure}

\begin{figure}[!htb]
\begin{tabular}{ccc}
\small $d\rho$ uniform measure  & \small  $d\rho$ Gaussian measure & \small $d\rho$ Chebyshev measure \\
\includegraphics[clip,width=\scalanellafigura\textwidth,viewport=0 0 1080 962]{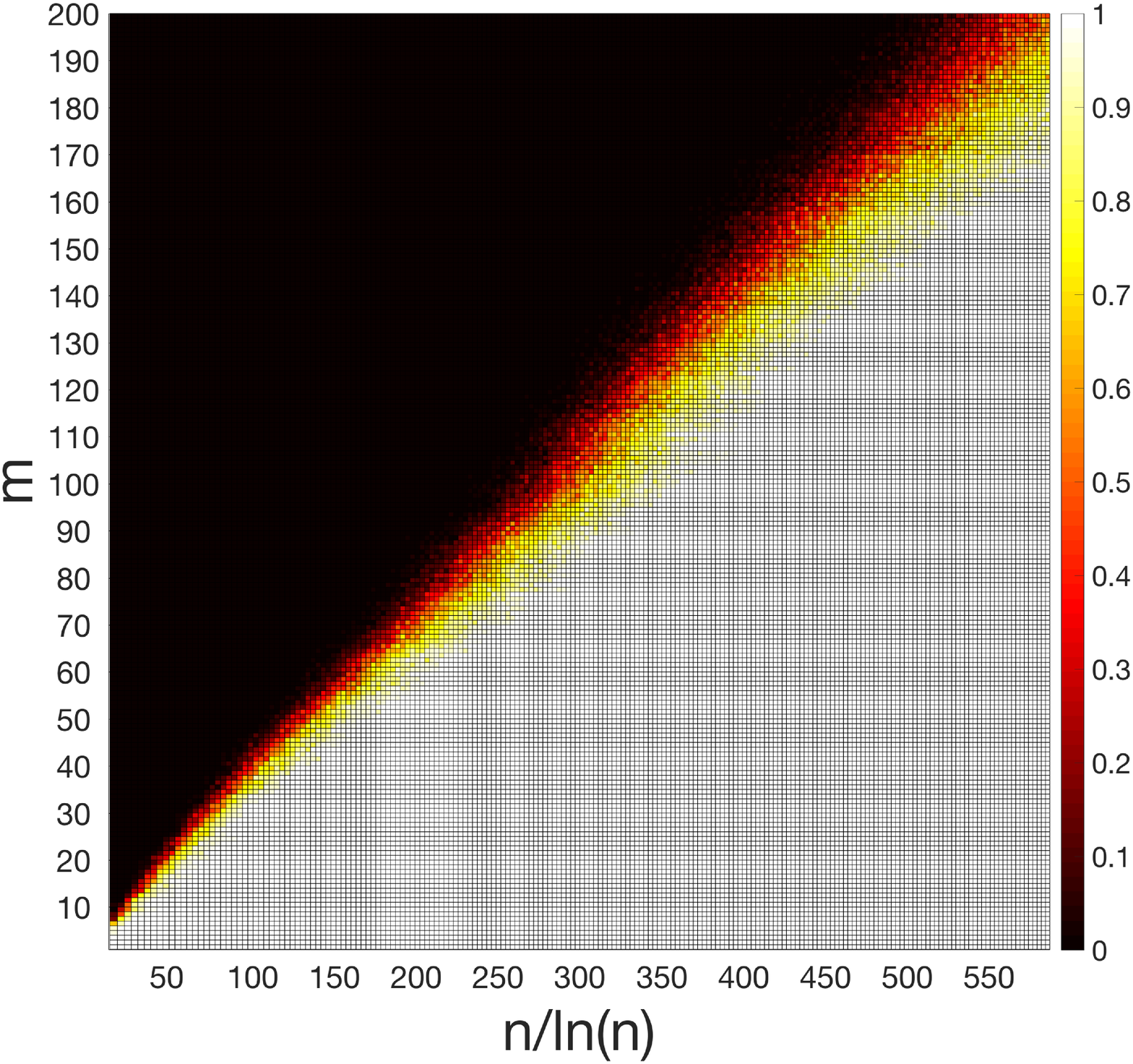}  &
\includegraphics[clip,width=\scalanellafigura\textwidth,viewport=0 0 1080 962]{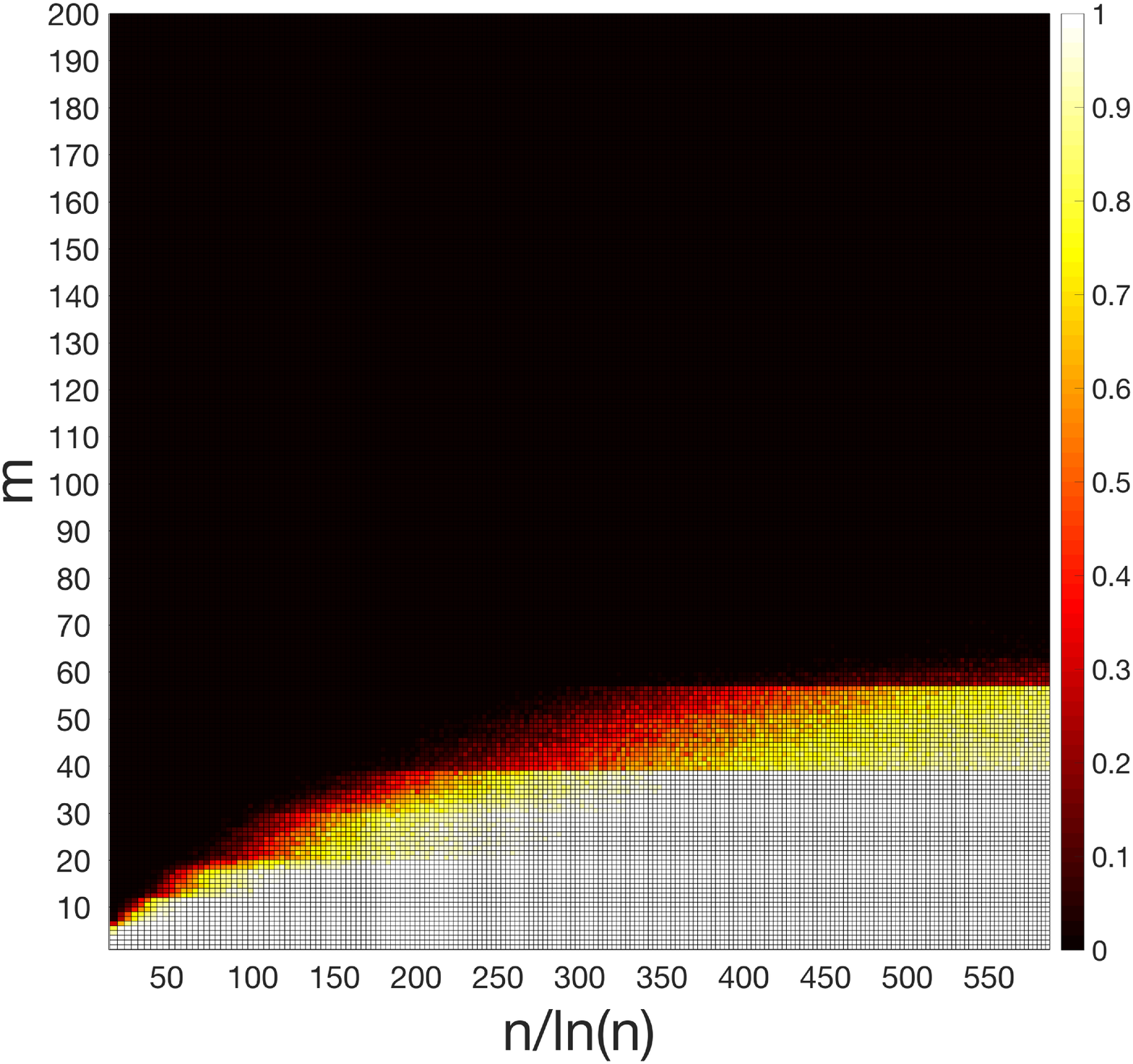}  &
\includegraphics[clip,width=\scalanellafigura\textwidth,viewport=0 0 1080 962]{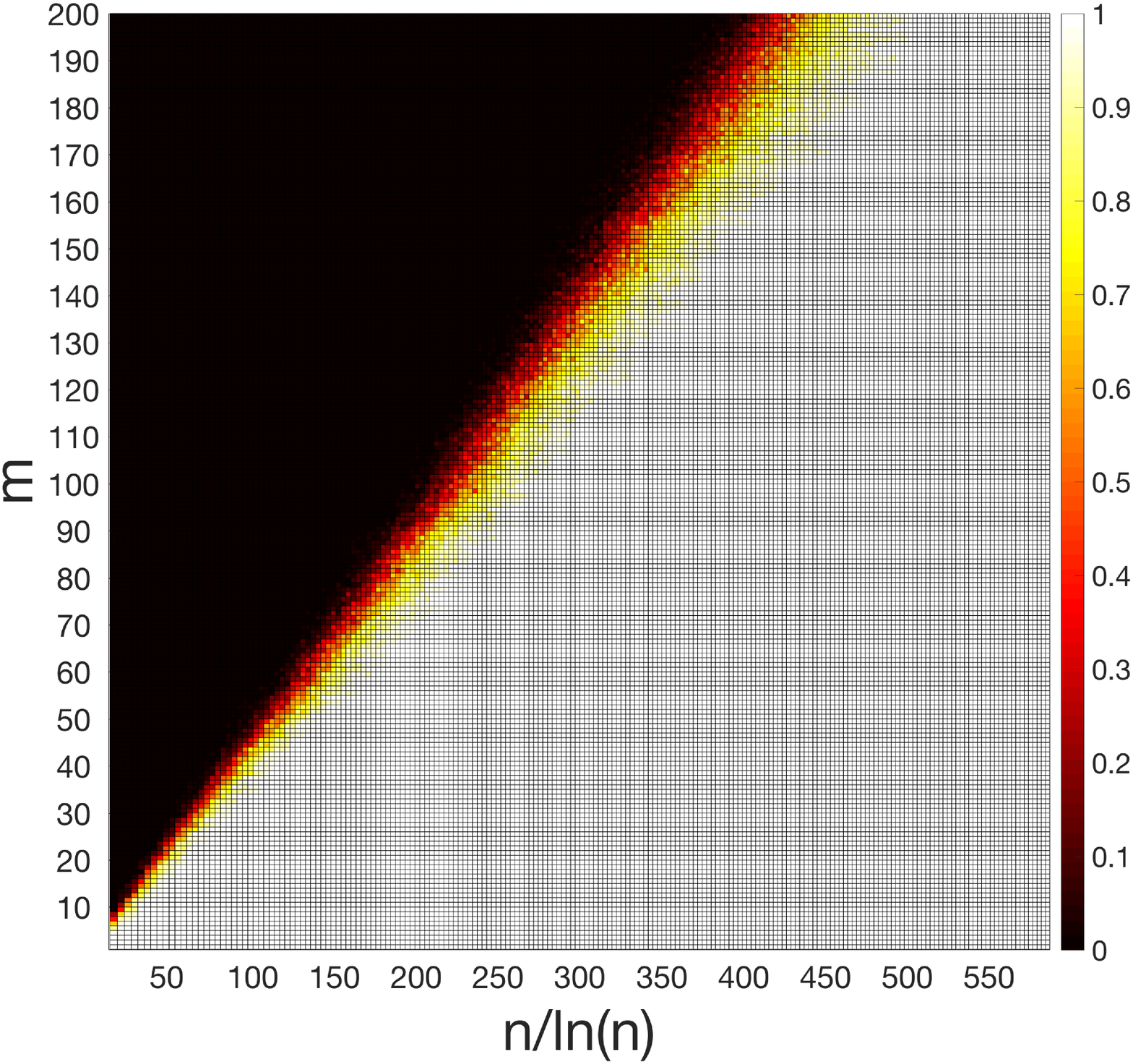}  
\end{tabular}
\caption{
Standard least squares, 
$Pr\{\textrm{cond}(G)\leq 3\}$, 
$d=10$. 
Left: $d\rho$ uniform measure. 
Center: $d\rho$ Gaussian measure.
Right: $d\rho$ Chebyshev measure. 
}
\label{fig_results_multivariate_d10_sls}
\end{figure}

Afterwards we present some numerical tests in the multivariate setting. 
These tests are again based, as in the previous section, on approximating the probability \eqref{eq:numerical_test_cond_prob} by empirical probability. 
In dimension $d$ larger than one there are many possible ways to enrich the polynomial space $\P_\Lambda$.  
The number of different downward closed sets whose cardinality equals $m$ gets very large already for moderate values of $m$ and $d$. 
Therefore, we present the numerical results for a chosen sequence of polynomial spaces $\P_{\Lambda_1},\ldots,\P_{\Lambda_m}$ 
such that  $\Lambda_1\subset \cdots \subset \Lambda_m$, where each $\Lambda_j\subset \N_0^d$ is downward closed, $\#(\Lambda_j)=\textrm{dim}(\P_{\Lambda_j})=j$ and the starting set $\Lambda_1$ contains only the null multi-index. 
All the tests in Fig.~\ref{fig_results_multivariate_d10_wls} and Fig.~\ref{fig_results_multivariate_d10_sls} have been obtained using the same sequence of increasingly embedded polynomial spaces $\P_{\Lambda_1}\subset \ldots \subset \P_{\Lambda_m}$, for both weighted and standard least squares and for the three choices of the measures $d\rho$. Such a choice allows us to establish a fair comparison between the two methods and among different measures, without the additional variability arising from modifications to the polynomial space. 

We report the results obtained for the tests in dimension $d=10$. 
The results in Fig.~\ref{fig_results_multivariate_d10_wls} confirm  
that weighted least squares always yield an empirical probability equal to one that $\textrm{cond}(\bG)\leq 3$, provided that $n/\log(n) \geq 2 m$.  
This condition ensures that \eqref{condmwc} with any choice of $r>0$ implies \eqref{eq:numerical_test_cond_prob}, thus verifying Corollary~\ref{cor1}. 
Again, the results do not show significant differences among the three choices of the measure $d\rho$: 
a straight line, with the same slope for all the three cases uniform, Chebyshev and Gaussian, separates the two regimes corresponding to empirical probabilities equal to zero and one.  
Compared to the univariate case in Fig.~\ref{fig_results_univariate_wls}, the results in Fig.~\ref{fig_results_multivariate_d10_wls} exhibit a sharper transition between the two extreme regimes, 
and an overall lower variability in the transition regime.

The results for standard least squares with $d=10$ are shown in Fig.~\ref{fig_results_multivariate_d10_sls}. In the case of the uniform measure, 
in Fig.~\ref{fig_results_multivariate_d10_sls}-right,
stability is ensured if $n/\ln(n) \geq 3.5 m$, 
which is more demanding than 
the condition $n/\ln(n) \geq 2 m$ needed for the stability of weighted least squares in Fig.~\ref{fig_results_multivariate_d10_wls}-right,  
but much less strict than 
the condition required 
with standard least squares 
in the univariate case, 
where $n/\ln(n)$ scales like $m^2$.  
These phenomena have already been observed and described in \cite{MNST}.
Similar results as those with the uniform measure are obtained with the Chebyshev measure in Fig.~\ref{fig_results_multivariate_d10_sls}-left, 
where again  
standard least squares achieve stability using more evaluations than weighted least squares in Fig.~\ref{fig_results_multivariate_d10_wls}-left.
The case of the Gaussian measure drastically differs from the uniform and Chebyshev cases: the results in Fig.~\ref{fig_results_multivariate_d10_sls}-center  clearly indicate that
a very large number of evaluations $n$ compared to $m$
is required to achieve   
stability of standard least squares. 

Let us mention that analogous results as those presented in Figs.~\ref{fig_results_univariate_wls} and \ref{fig_results_multivariate_d10_wls}
for weighted least squares have been obtained also in other dimensions, 
and with many other sequences of increasingly embedded polynomial spaces.  
In the next tables we report some of these results for selected values of $d=1,2,5,10,50,100$. 
We choose $n=26599$ and $m=200$ that satisfy condition \eqref{condmwc} with $r=1$, 
and report in Table~\ref{tab:results_high_dimensions} the empirical probabilities that approximate \eqref{eq:numerical_test_cond_prob}, 
again calculated over one hundred repetitions. 
This table provides multiple comparisons: weighted least squares versus standard least squares, for the three choices of the measure $d\rho$ (uniform, Gaussian and Chebyshev) and with $d$ varying between $1$ and $100$.

\begin{table}[!htb]
\centering
\begin{tabular}{|c|c|c|c|c|c|c|c|c|}
  \hline
 method     &   $d\rho$  &  $d=1$  & $d=2$  & $d=5$ & $d=10$   & $d=50$ & $d=100$ \\
  \hline
\hline
\hline
weighted LS & uniform    &    1     &   1    &    1  &    1    &   1    &    1    \\
weighted LS & Gaussian   &    1     &   1    &    1  &    1    &   1    &    1    \\
weighted LS & Chebyshev  &    1     &   1    &    1  &    1    &   1    &    1    \\
\hline
\hline
standard LS & uniform    &    0     &   0    &    0.54  &    1    &   1    &    1    \\
standard LS & Gaussian   &    0     &   0    &    0  &    0    &   0    &    0    \\
standard LS & Chebyshev  &    1     &   1    &    1  &    1    &   1    &    1    \\   
  \hline
\end{tabular}
\caption{$Pr\{\textrm{cond}(G)\leq 3\}$, with $n=26559$ and $m=200$: weighted least squares versus standard least squares, $d\rho$ uniform versus $d\rho$ Gaussian versus $d\rho$ Chebyshev, $d=1, 2, 5, 10, 50, 100$.
}
\label{tab:results_high_dimensions}
\end{table}
\begin{table}[!htb]
\centering
\begin{tabular}{|c|c|c|c|c|c|c|c|c|}
  \hline
 method     &   $d\rho$  &  $d=1$  & $d=2$  & $d=5$ & $d=10$   & $d=50$ & $d=100$ \\
  \hline
\hline
\hline
weighted LS & uniform   & $1.5593$  & $1.4989$  & $1.4407$ & $1.4320$ & $1.4535$ & $1.4179$   \\
weighted LS & Gaussian  & $1.5994$  & $1.5698$  & $1.4743$ & $1.4643$ & $1.4676$ & $1.4237$   \\
weighted LS & Chebyshev & $1.5364$  & $1.4894$  & $1.4694$ & $1.4105$ & $1.4143$ & $1.4216$   \\
\hline
\hline
standard LS & uniform   & $19.9584$ & $29.8920$ & $3.0847$ & $1.9555$ & $1.7228$ & $1.5862$  \\
standard LS & Gaussian  & $\sim10^{19}$ & $\sim 10^{19}$  & $ \sim 10^{19}$ & $\sim 10^{16}$ & $ \sim 10^{9}$ &$ \sim 10^{3}$\\
standard LS & Chebyshev & $1.5574$ & $1.5367$ & $1.5357$ & $1.4752$ & $1.4499$ & $1.4625$ \\   
  \hline
\end{tabular}
\caption{Average of $\textrm{cond}(G)$, with $n=26559$ and $m=200$: weighted least squares versus standard least squares, $d\rho$ uniform versus $d\rho$ Gaussian versus $d\rho$ Chebyshev, $d=1, 2, 5, 10, 50, 100$.
}
\label{tab:results_high_dimensions_cond}
\end{table}

In Table~\ref{tab:results_high_dimensions}, all the empirical probabilities related to results for weighted least squares are equal to one, and confirm the theory since, for the chosen values of $n$, $m$ and $r$, the probability \eqref{eq:numerical_test_cond_prob}  is larger than $1-5.67 \times 10^{-7}$. 
This value is computed using estimate \eqref{eq:finer_est_prob} from the proof of Theorem~\ref{theo2}.
In contrast to weighted least squares, whose empirical probability equal one independently of $d\rho$ and $d$, 
the empirical probability of standard least squares does depend on the chosen measure, and to some extent on the dimension $d$ as well.   
With the uniform measure, the empirical probability that approximates \eqref{eq:numerical_test_cond_prob} equals zero when $d=1$ or $d=2$, equals $0.54$ when $d=5$, and equals one when $d=10$, $d=50$ or $d=100$. 
In the Gaussian case, standard least squares always feature null empirical probabilities. 
With the Chebyshev measure, the condition number of $\bG$ for standard least squares is always lower than three for any tested value of $d$.

In addition to the results in Table~\ref{tab:results_high_dimensions}, 
further information 
are needed for assessing how severe is the lack of stability when obtaining null empirical probabilities.  
To this aim, in Table~\ref{tab:results_high_dimensions_cond} we also report the average value of $\textrm{cond}(\bG)$, obtained when averaging the condition number of $\bG$ over the same repetitions used to estimate the empirical probabilities in Table~\ref{tab:results_high_dimensions}. 
The information in Table~\ref{tab:results_high_dimensions_cond} are complementary to those in Table~\ref{tab:results_high_dimensions}. 
On the one hand they point out the stability and robustness of weighted least squares, showing a tamed condition number with any measure $d\rho$ and any dimension $d$.   
On the other hand they provide further insights on stability issues of standard least squares and their dependence on $d\rho$ and $d$. 
For standard least squares with the uniform measure, the average condition number reduces as the dimension $d$ increases, in agreement with the conclusion drawn from Table~\ref{tab:results_high_dimensions}. 
The Gramian matrix of standard least squares with the Gaussian measure is very ill-conditioned for all tested values of $d$. 
For standard least squares with the Chebyshev measure, the averaged condition number of $\bG$ is only slightly larger than the one for weighted least squares. 

It is worth remarking that, the results for standard least squares in Fig.~\ref{fig_results_multivariate_d10_sls}, Table~\ref{tab:results_high_dimensions} and Table~\ref{tab:results_high_dimensions_cond} are sensitive to the chosen sequence of polynomial spaces.
Testing different sequences might produce different results, that however necessarily obey to the estimates proven in Theorem~\ref{theo1} with uniform and Chebyshev measures, when $n$, $m$ and $r$ satisfy condition \eqref{condm}.
Many other examples with 
standard least squares have been extensively discussed in previous works \emph{e.g.}~\cite{MNST,CCMNT}, 
also in situations where $n$, $m$ and $r$ do not satisfy condition \eqref{condm} and therefore Theorem~\ref{theo1} does not apply.
In general, when $n$, $m$ and $r$ do not satisfy \eqref{condm}, 
there exist multivariate polynomial spaces 
of dimension $m$ such that the Gramian matrix of standard least squares with the uniform and Chebyshev measures does not satisfy \eqref{tailhalf}. 
Examples of such spaces are discussed in \cite{MNST,CCMNT}. 
Using these spaces would yield null empirical probabilities in Table~\ref{tab:results_high_dimensions} for standard least squares with the uniform and Chebyshev measures.

For weighted least squares, when $n$, $m$ and $r$ satisfy condition \eqref{condmwc}, 
any sequence of polynomial spaces 
yields empirical probabilities close to one, 
according to Corollary~\ref{cor1}. 
Indeed such a robustness with respect to the choices of $d\rho$, of the polynomial space and of the dimension $d$ represents one of the main advantages of the weighted approach.

\end{document}